\documentclass[11pt,a4paper]{article}

 \usepackage{amssymb,amsmath,amsfonts,amssymb}
 \usepackage{graphics,graphicx,color}
 -\marginparwidth \oddsidemargin -4mm \evensidemargin
 \oddsidemargin%
 \advance\hoffset by 5mm
 \usepackage{amsmath}    
 \advance\hoffset by 5mm

\newcommand\grad{{\bf \nabla}}

\newcommand\la{{\lambda}}
\newcommand\ga{{\gamma}}
\newcommand\nvec{{\bf n}}
\newcommand\xvec{{\bf x}}
\newcommand\uvec{{\bf u}}
\newcommand\yvec{{\bf y}}
\newcommand\evec{{\bf e}}
\newcommand\fvec{{\bf f}}
\newcommand\Bvec{{\bf B}}
\newcommand\Svec{{\bf S}}

\newcommand\mvec{{\bf m}}

\newcommand\Qvec{{\bf Q}}
\newcommand\Hvec{{\bf H}}
\newcommand \I {{\bf I_{LG}}}
\newcommand \Ihat {{\bf \hat{I}_{LG}}}

\newcommand{\Acal}{{\cal A}}

\newcommand{\Rr}{{\mathbb R}}

 \title{Symmetry of uniaxial global Landau-de Gennes minimizers in the theory of nematic liquid crystals}
 \date{\today}
\author{Duvan Henao and Apala Majumdar}

 \newtheorem{lemma} {Lemma}
 \newtheorem{proposition}{Proposition}
 \newtheorem{theorem} {Theorem}

\begin{document}

\maketitle


\section{Preliminaries}
\label{sec:prelim} Let $B\left(0,R_0\right)\subset \Rr^3$ denote a
three-dimensional spherical droplet of radius $R_0>0$, centered at
the origin. Let $S_0$ denote the set of symmetric, traceless
$3\times 3$ matrices i.e.
\begin{equation}
\label{eq:1} S_0 = \left\{\Qvec\in M^{3\times 3};
\Qvec_{ij}=\Qvec_{ji}; \Qvec_{ii}=0 \right\}
\end{equation}
where $M^{3\times 3}$ is the set of $3\times 3$ matrices. The
corresponding matrix norm is defined to be \cite{amaz}
\begin{equation}
\label{eq:2} |\Qvec|^2 = \Qvec_{ij} \Qvec_{ij} \quad i,j=1\ldots 3
\end{equation}
and we will use the Einstein summation convention throughout the
paper.

We work with the Landau-de Gennes theory for nematic liquid
crystals \cite{dg} whereby a liquid crystal configuration is
described by a macroscopic order parameter, known as the
$\Qvec$-tensor order parameter. Mathematically, the Landau-de
Gennes $\Qvec$-tensor order parameter is a symmetric, traceless
$3\times 3$ matrix belonging to the space $S_0$ in (\ref{eq:1}).
The liquid crystal energy is given by the Landau-de Gennes energy
functional and the associated energy density is a nonlinear
function of $\Qvec$ and its spatial derivatives
\cite{dg,newtonmottram}. We work with the simplest form of the
Landau-de Gennes energy functional that allows for a first-order
nematic-isotropic phase transition and spatial inhomogeneities as
shown below \cite{amaz} -
\begin{equation}
\label{eq:3} \I \left[\Qvec \right] = \int_{B(0,R_0)}
\frac{L}{2}|\grad \Qvec|^2 + f_B\left(\Qvec\right)~dV.
\end{equation}
Here, $L>0$ is a small material-dependent elastic constant,
$|\grad \Qvec|^2 = \Qvec_{ij,k}\Qvec_{ij,k}$ ( note that $\Qvec_{ij,k} = \frac{\partial \Qvec_{ij}}{\partial \xvec_k}$) with $i,j,k=1\ldots
3$ is an \emph{elastic energy density} and $f_B:S_0 \to \Rr$ is
the \emph{bulk energy density}. For our purposes, we take $f_B$ to
be a quartic polynomial in the $\Qvec$-tensor invariants as shown
below -
\begin{equation}
\label{eq:4} f_B(\Qvec) = \frac{\alpha (T -
T^*)}{2}\textrm{tr}\Qvec^2 - \frac{b^2}{3}\textrm{tr}\Qvec^3 +
\frac{c^2}{4}\left(\textrm{tr}\Qvec^2\right)^2
\end{equation}
where $\textrm{tr}\Qvec^3 = \Qvec_{ij}\Qvec_{jp}\Qvec_{pi}$ with
$i,j,p=1\ldots 3$, $\alpha, b^2, c^2>0$ are material-dependent
constants, $T$ is the absolute temperature and $T^*$ is a
characteristic temperature below which the isotropic phase
$\Qvec=0$ loses its stability. We work in the low-temperature
regime with $T<<T^*$ and hence, we can re-write (\ref{eq:4}) as
\begin{equation}
\label{eq:5} f_B(\Qvec) = -\frac{a^2}{2}\textrm{tr}\Qvec^2 -
\frac{b^2}{3}\textrm{tr}\Qvec^3 +
\frac{c^2}{4}\left(\textrm{tr}\Qvec^2\right)^2
\end{equation} where $a^2>0$ is a temperature-dependent parameter
and we will subsequently investigate the $a^2\to \infty$ limit,
known as the \emph{low-temperature} limit. One can readily verify
that $f_B$ is bounded from below and attains its minimum on the
set of $\Qvec$-tensors given by \cite{ejam,am2}
\begin{equation}
\label{eq:6} \Qvec_{\min} = \left\{ \Qvec \in S_0; \Qvec = s_+
\left(\nvec \otimes \nvec - \frac{\mathbf{I}}{3}\right),~\nvec \in
S^2 \right\}
\end{equation}
where $\mathbf{I}$ is the $3\times 3$ identity matrix and 
\begin{equation}
\label{eq:s+} s_+ = \frac{b^2 + \sqrt{b^4 + 24 a^2 c^2}}{4 c^2}.
\end{equation}

We are interested in characterizing global minimizers of the
Landau-de Gennes energy functional in (\ref{eq:3}), on spherical
droplets with \emph{homeotropic} or \emph{radial anchoring
conditions} \cite{am2}. The global Landau-de Gennes minimizers
correspond to physically observable liquid crystal configurations
and hence, are of both mathematical and practical importance. We
take our admissible $\Qvec$-tensors to belong to the space
\begin{equation}
\label{eq:7} \Acal = \left\{ \Qvec \in
W^{1,2}\left(B(0,R_0);S_0\right); \Qvec = \Qvec_b ~on~
\partial B(0,R_0) \right\}
\end{equation}
where $W^{1,2}\left(B(0,R_0);S_0\right)$ is the Soboblev space of
square-integrable $\Qvec$-tensors with square-integrable first
derivatives \cite{evans}. The Dirichlet boundary condition
$\Qvec_b$ is given by
\begin{equation}
\label{eq:8} \Qvec_b(\xvec) = s_+ \left( \frac{\xvec \otimes
\xvec}{|\xvec|^2} - \frac{\mathbf{I}}{3}\right) \in \Qvec_{\min}
\end{equation}
where $\xvec \in \Rr^3$ is the position vector and
$\frac{\xvec}{|\xvec|}$ is the unit-vector in the radial
direction. For completeness, we recall that the $W^{1,2}$-norm is
given by $||\Qvec||_{W^{1,2}} = \left(\int_{B(0,R_0)} |\Qvec|^2 +
|\grad \Qvec|^2~dV \right)^{1/2}$ and the $L^{\infty}$-norm is
defined to be $||\Qvec||_{L^{\infty}} = \textrm{ess sup}_{\xvec
\in B(0,R_0)}|\Qvec (\xvec)|$ \cite{evans}.

We define a modified Landau-de Gennes energy functional as shown
below -
\begin{equation}
\label{eq:9} \Ihat \left[\Qvec \right] = \int_{B(0,R_0)}
\frac{L}{2}|\grad \Qvec|^2 + f_B\left(\Qvec\right) -
\min_{\Qvec\in S_0}f_B\left(\Qvec\right) ~dV
 \end{equation}
 where $f_B\left(\Qvec\right) -
\min_{\Qvec\in S_0}f_B\left(\Qvec\right) \geq 0$ for all $\Qvec\in
S_0$. The existence of a global minimizer for the modified
functional $\Ihat$ in the admissible space $\Acal$ is immediate
from the direct methods in the calculus of variations
\cite{evans}; the details are omitted for brevity. It is clear
that $\Qvec^*\in\Acal$ is a minimizer of $\Ihat$ if and only if
$\Qvec^*$ is a minimizer of $\I$ in (\ref{eq:3}) and hence, it
suffices to study minimizers of the modified functional in
(\ref{eq:9}). In what follows, we study the modified functional in
(\ref{eq:9}) and drop the \emph{hat} for brevity.

For a fixed $a^2>0$, let $\Qvec^a \in \Acal$ denote a
corresponding Landau-de Gennes global minimizer. The
Euler-Lagrange equations associated with $\Ihat$ are given by a
nonlinear elliptic system of coupled partial differential
equations:
\begin{equation}
\label{eq:10} L \Delta \Qvec_{ij} = -a^2 \Qvec_{ij} -
b^2\left(\Qvec_{ip}\Qvec_{pj} -
\frac{1}{3}\Qvec_{pq}\Qvec_{pq}\delta_{ij}\right) + c^2
\left(\textrm{tr}\Qvec^2\right) \Qvec_{ij} \quad i,j,p,q=1\ldots 3
\end{equation} where
$\frac{b^2}{3}\Qvec_{pq}\Qvec_{pq}\delta_{ij}$ is a Lagrange
multiplier accounting for the tracelessness constraint
\cite{amaz}. It follows from standard arguments in elliptic
regularity that $\Qvec^a$ is smooth and real analytic on
$B(0,R_0)$.

The notion of a \emph{limiting harmonic map} was first introduced
in \cite{amaz} and is crucial in what follows. A limiting harmonic
map $\Qvec^0 \in \Acal$ is defined to be
\begin{equation}
\label{eq:11} \Qvec^0 = s_+ \left( \nvec^0 \otimes \nvec^0 -
\frac{\mathbf{I}}{3} \right)
\end{equation}
where $\nvec^0$ is a minimizer of the Dirichlet energy
\cite{schoen}
\begin{equation}
\label{eq:12} I[\nvec] = \int_{B(0,R_0)} |\grad \nvec|^2~dV
\end{equation}
in the admissible space $\Acal_\nvec = \left\{ \nvec \in
W^{1,2}\left(B(0,R_0);S^2\right); \nvec = \frac{\xvec}{|\xvec|}
~on~ \partial B(0,R_0)\right\}$. In the case of a spherical
droplet with homeotropic boundary conditions, $\nvec^0$ is unique
and given by the radial unit-vector \cite{lin}. Hence, the
limiting harmonic map is unique for our model problem and is given by
\begin{equation}
\label{eq:13} \Qvec^0 = s_+ \left(\frac{\xvec\otimes
\xvec}{|\xvec|^2} - \frac{\mathbf{I}}{3} \right).
\end{equation} We note that $\Qvec^0$ has a single isolated point
defect at the origin.

In what follows, we keep $L, b^2$ and $c^2$ fixed in (\ref{eq:9}).
 Following the methods in \cite{am2}, we
first introduce a re-scaling relevant to the
\emph{low-temperature} limit: $a^2\to\infty$. Let $t$ denote a
\emph{dimensionless} temperature
\begin{equation}
\label{eq:t} t = \frac{27 a^2 c^2}{b^4}>0
\end{equation}
so that the $a^2\to\infty$ limit corresponds to the $t\to\infty$
limit and define
\begin{equation}
\label{eq:h+} h_+ = \frac{3 + \sqrt{9 + 8t}}{4} \sim
\sqrt{\frac{t}{2}} \quad t\to\infty.
\end{equation} Define
\begin{equation}
\label{eq:14} \bar{\Qvec}_{ij} =
\frac{1}{h_+}\sqrt{\frac{27 c^4}{2 b^4}}\Qvec_{ij}
\end{equation}
and the corresponding Landau-de Gennes energy functional is given
by (up to a multiplicative constant)
\begin{equation}
\label{eq:15} \Ihat[\bar{\Qvec}] = \int_{B(0,R_0)}
\frac{\bar{L}}{2}|\grad \bar{\Qvec}|^2 + \frac{t}{8}\left[
\left(1- |\bar{\Qvec}|^2 \right)^2 - \frac{8 h_+}{t}
\sqrt{\frac{3}{2}}\textrm{tr}\bar{\Qvec}^3+ f(t)\right]~dV
\end{equation}
where $\bar{L} = \frac{27 c^2 L}{2 b^4}>0$ is fixed, $f(t)$ is a
function that can be explicitly computed, \newline  $\left(1-
|\bar{\Qvec}|^2 \right)^2 - \frac{8 h_+}{t}
\sqrt{\frac{3}{2}}\textrm{tr}\bar{\Qvec}^3+ f(t) \geq 0$ for
$\Qvec\in S_0$ (see definition of $\Ihat$ in (\ref{eq:9})) and for
$t$ sufficiently large, we have
\begin{equation}
\label{eq:16} \frac{\sigma_1}{\sqrt{t}}\leq \frac{h_+}{t} \leq
\frac{\sigma_2}{\sqrt{t}} ;\quad \frac{\gamma_1}{\sqrt{t}}\leq
f(t) \leq \frac{\gamma_2}{\sqrt{t}} \quad t\to\infty
\end{equation}
for positive constants $\sigma_1, \sigma_2$ and constants
$\gamma_1, \gamma_2$ independent of $t$ in the $t\to\infty$ limit.
The corresponding admissible space is
\begin{equation}
\label{eq:7new} \bar{\Acal} = \left\{ \bar{\Qvec} \in
W^{1,2}\left(B(0,R_0);S_0\right); \bar{\Qvec} =
\frac{1}{h_+}\sqrt{\frac{27 c^4}{2 b^4}}\Qvec_b ~on~
\partial B(0,R_0) \right\}
\end{equation} where $\Qvec_b$ has been defined above. In what follows, we drop the
\emph{bars} from the re-scaled variables for brevity and all
subsequent statements in this section are to be understood in
terms of the scaled variable in (\ref{eq:14}).

 Next, we
quote some results from \cite{amaz} that are relevant to
the development of our mathematical framework.

\begin{proposition}
\label{prop:1} For each $t>0$, let $\Qvec^t=s^t\left(\nvec^t \otimes \nvec^t - \frac{1}{3}\mathbf{I} \right)$ denote a uniaxial global
minimizer of $\Ihat$ in (\ref{eq:15}), in the admissible space
$\bar{\Acal}$ defined in (\ref{eq:7new}), where $s^t \in W^{1,2}\left(\Omega, \Rr\right)$ and $\nvec^t\in W^{1,2}\left(\Omega, S^2 \right)$. Then there exists a
sequence $\left\{t_k\right\}$ with $t_k\to\infty$ as $k\to\infty$
such that $\Qvec^{t_k} \to \Qvec^0$ strongly in $W^{1,2}\left(
B(0,R_0);S_0 \right)$, where $\Qvec^0$ is the re-scaled limiting harmonic
map defined in (\ref{eq:13}).
\end{proposition}

\textit{Proof:} The proof of Proposition~\ref{prop:1} is identical
to the proof of Lemma~3 in \cite{amaz} and the details are omitted
for brevity. From the global energy minimality, we have
\begin{equation}
\label{eq:18} \Ihat[\Qvec^{t_k}] \leq \Ihat[\Qvec^0] = E\left(L,
R_0\right)~ \forall k
\end{equation} where $E>0$ is a positive constant independent of $t_k$ i.e. the $\Qvec^{t_k}$'s
have bounded energy in the limit $k\to\infty$. Note that
$$ \frac{t}{8}\left[
\left(1- |\Qvec^0|^2 \right)^2 - \frac{8 h_+}{t}
\sqrt{\frac{3}{2}}\textrm{tr}\left(\Qvec^0\right)^3+ f(t)\right] =
0$$ for all $t>0$, by definition (since $\Qvec^0 \in \Qvec_{\min}$). One immediate consequence of the
strong convergence result are the following integral equalities -
\begin{eqnarray}
\label{eq:17} && \int_{B(0,R_0)}|\grad \Qvec^{0}|^2~dV =
\lim_{k\to\infty} \int_{B(0,R_0)}|\grad \Qvec^{t_k}|^2~dV
\nonumber \\ && \lim_{k\to\infty} \int_{B(0,R_0)}\left[ \left(1-
|\Qvec^{t_k}| \right)^2 - \frac{8 h_+}{t}
\sqrt{\frac{3}{2}}\textrm{tr}\left(\Qvec^{t_k}\right)^3+
f(t)\right]~dV =0.
\end{eqnarray} From the inequalities (\ref{eq:16}), we deduce that
\begin{equation}
\label{eq:19}\lim_{k\to\infty}\int_{B(0,R_0)} \left(1-
|\Qvec^{t_k}|^2 \right)^2~dV = 0. \end{equation} $\Box$

\begin{proposition}
\label{prop:3}For each $t>0$, let $\Qvec^t$ denote a global
minimizer of $\Ihat$ in (\ref{eq:15}), in the admissible space
$\bar{\Acal}$ defined in (\ref{eq:7new}). Let $\left\{t^k\right\}$
be a sequence such that $t^k \to \infty$ as $k\to\infty$. Then
\begin{equation}
\label{eq:21} \lim_{k\to\infty} |\Qvec^{t_k}(\xvec)| \leq 1 \quad
\forall \xvec \in B(0,R_0).
\end{equation}
\end{proposition}
\textit{Proof:} The proof follows from a standard maximum
principle; see \cite{amaz} Proposition 3 for an analogous
statement with proof. $\Box$

\begin{proposition}
\label{prop:2} [\cite{amaz}; ~ Lemma $2$] For each $t>0$, let $\Qvec^t$
denote a global minimizer of $\Ihat$ in (\ref{eq:15}), in the
admissible space $\bar{\Acal}$ defined in (\ref{eq:7new}). Define
$B(\xvec,r) = \left\{ \yvec \in B(0,R_0); \left| \xvec - \yvec
\right| \leq r \right\} \subset B(0,R_0)$ and $$
e\left(\Qvec^t,\grad \Qvec^t\right) = \frac{\bar{L}}{2}|\grad
\Qvec^t|^2 + \frac{t}{8}\left[ \left(1- |\Qvec^t|^2 \right)^2 -
\frac{8 h_+}{t}
\sqrt{\frac{3}{2}}\textrm{tr}\left(\Qvec^t\right)^3+
f(t)\right].$$ Then
\begin{equation}
\label{eq:20} \frac{1}{r}\int_{B(\xvec,r)} e\left(\Qvec^t,\grad
\Qvec^t\right)~ dV \leq \frac{1}{R}\int_{B(\xvec,R)}
e\left(\Qvec^t,\grad \Qvec^t\right)~ dV, ~ \forall \xvec \in
B(0,R_0),~ r\leq R,
\end{equation} so that $B(\xvec, R)\subset B(0,
R_0)$.
\end{proposition}
\textit{Proof:} The proof can be found in [\cite{amaz};~ Lemma $2$].
An analogous boundary monotonicity formula can be found in Lemma 9
\cite{amaz}. $\Box$

\begin{proposition}
\label{prop:4}For each $t>0$, let $\Qvec^t$ denote a global uniaxial
minimizer of $\Ihat$ in (\ref{eq:15}), in the admissible space
$\bar{\Acal}$ defined in (\ref{eq:7new}). Let $\left\{t^k\right\}$
be a sequence such that $t^k \to \infty$ as $k\to\infty$ with
$\Qvec^{t_k} \to \Qvec^0$ in $W^{1,2}\left(B(0,R_0),S_0\right)$ as
$k\to\infty$, where $\Qvec^0$ is the re-scaled limiting harmonic map defined
in (\ref{eq:13}). For any compact $K\subset B(0, R_0)$ such that
$K$ does not contain any singularities of $\Qvec^0$ i.e. does not
contain the origin, we have
\begin{equation}
\label{eq:22} \lim_{k\to\infty}|\Qvec^{t_k}(\xvec)| = 1 \quad
\forall \xvec \in K
\end{equation}
and the limit is uniform on $K$.
\end{proposition}
\textit{Proof:} Proposition~\ref{prop:4} is a consequence of the
pointwise uniform convergence
$$ \lim_{k\to\infty}\left[ \left(1-
|\Qvec^{t_k}| \right)^2 - \frac{8 h_+}{t}
\sqrt{\frac{3}{2}}\textrm{tr}\left(\Qvec^{t_k}\right)^3+
f(t_k)\right] = 0$$ everywhere away from the singular set of
$\Qvec^0$ i.e. away from the origin. This uniform convergence result holds in the interior
and up to the boundary. The proof can be found in \cite{amaz},
Propositions $4$ and $6$. $\Box$

To summarize, let $\left\{t^k\right\}$ be a sequence such that
$t^k \to \infty$ as $k\to\infty$ and let $\left\{\Qvec^{t_k} \right\}$ denote a corresponding sequence of purely uniaxial Landau-de Gennes minimizers. Then up to a subsequence,
$\Qvec^{t_k} \to \Qvec^0$ in $W^{1,2}\left(B(0,R_0),S_0\right)$ as
$k\to\infty$. Further for $k$ sufficiently large,
$\Ihat[\Qvec^{t_k}]$ can be bounded independently of $t_k$ since
$$\Ihat[\Qvec^{t_k}] \leq \Ihat[\Qvec^0]$$
and $|\Qvec^{t_k}|$ is strictly positive (and bounded from below)
everywhere away from the origin.

\section{Statement of main results}
\label{sec:main}

Our main result is the following :
\begin{theorem}
\label{thm:1} Let $B(0,R_0)\subset \Rr^3$ denote a spherical
droplet of radius, $R_0$, centered at the origin. For each $a^2>0$,
let $\Qvec^a$ denote a global Landau-de Gennes minimizer in the
admissible space $\Acal$ defined in (\ref{eq:7}). Then in the
limit $a^2\to\infty$, $\Qvec^a$ cannot be purely uniaxial i.e.
cannot be of the form
\begin{equation}
\label{eq:23} \Qvec^a(\xvec) = s(\xvec) \left( \nvec \otimes \nvec
- \frac{\mathbf{I}}{3}\right) \quad \xvec\in B(0,R_0)
\end{equation}
for a function $s: B(0,R_0) \to \Rr$ and a unit-vector field
$\nvec\in W^{1,2}\left(B(0,R_0),S^2\right)$.
\end{theorem}

\textit{Comment: From \cite{bz}, it is known that if $\Qvec=s\left(\nvec\otimes\nvec - \frac{\mathbf{I}}{3} \right) \in W^{1,2}\left(\Omega, S_0 \right)$ on a simply-connected domain $\Omega \subset \Rr^3$, then $\nvec\in W^{1,2}\left(\Omega, S^2 \right)$.}

We prove Theorem~\ref{thm:1} by contradiction in the subsequent
sections. Of key importance is the division trick used in
\cite{pisante} for the Ginzburg-Landau theory for
superconductivity in three dimensions. We adapt the division trick
in \cite{pisante} to the Landau-de Gennes framework for
nematic liquid crystals. We point out the following important
differences: (i) we have a parameter $a^2$ and we are interested in
the asymptotics of global energy minimizers in the $a^2\to\infty$
limit; the Ginzburg-Landau equations in \cite{pisante} are parameter-free, (ii) the nonlinearities in the Landau-de
Gennes equations (\ref{eq:10}) are different to the nonlinearities
in the Ginzburg-Landau equations introducing additional technical
complexities and (iii) the Landau-de Gennes macroscopic variable is a two-tensor field
$\Qvec\in S_0$ whereas the Ginzburg-Landau macroscopic variable in
$\Rr^3$ is a three-dimensional vector field $\uvec \in \Rr^3$.

To make better contact with the framework used in \cite{pisante},
we introduce the following dimensionless variables as in
\cite{am2} :
\begin{eqnarray}
\label{eq:nondim} && \tilde{\xvec} = \frac{\xvec}{\xi_b}, \quad
\tilde{\Qvec} = \frac{1}{h_+}\sqrt{\frac{27 c^4}{ 2b^4}} \Qvec, ~
\tilde{\mathcal{I}}_{LG} = \frac{h_+^2}{\sqrt{t}} \sqrt{\left(\frac{27
c^6}{4 b^4 L^3}\right)}\Ihat
\end{eqnarray}
where $t = \frac{27 a^2 c^2}{b^4}>0$ is the \emph{reduced
temperature} \cite{mkaddem&gartland2}, $t>1$ throughout the paper
and $h_+$ has been defined in (\ref{eq:h+}). The length-scale
$\xi_b = \sqrt{\frac{27c^2 L}{ t b^4}}$. We note that the position
vector $\xvec$ has been re-scaled in (\ref{eq:nondim}) whereas it
was left unchanged in (\ref{eq:14}). The corresponding
dimensionless energy density is
\begin{eqnarray}
\label{eq:nondim3} \tilde{e}(\tilde{\Qvec},\grad\tilde{\Qvec}) =
\frac{1}{2}|\grad\tilde{\Qvec}|^2 -
\frac{1}{2}\textrm{tr}\tilde{\Qvec}^2 -
\frac{\sqrt{6}h_+}{t}\textrm{tr}\tilde{\Qvec}^3 + \frac{h_+^2}{2t}
\left(\textrm{tr}\tilde{\Qvec}^2\right)^2 + C(t)
\end{eqnarray}
where $C(t) = \frac{1}{2} + \frac{h_+}{t} - \frac{h_+^2}{2t}$ is
an additive constant that ensures
$$- \frac{1}{2}\textrm{tr}\tilde{\Qvec}^2 -
\frac{\sqrt{6}h_+}{t}\textrm{tr}\tilde{\Qvec}^3 + \frac{h_+^2}{2t}
\left(\textrm{tr}\tilde{\Qvec}^2\right)^2 + C(t) \geq 0 \quad \forall \Qvec \in S_0.$$ The
corresponding Landau-de Gennes energy functional is given by
\begin{equation}
\label{eq:nondim4} \tilde{\mathcal{I}}_{LG}[\tilde{\Qvec}] =
\int_{B(0,\tilde{R}_t)}\tilde{e}(\tilde{\Qvec},\grad\tilde{\Qvec})~dV,
\end{equation}
where $\tilde{R}_t = \gamma \sqrt{t}R_0$ for a fixed constant $\gamma>0$. In particular,
$\tilde{R}_t \to \infty$ as $t \to \infty$. In what follows, we
drop the
 \emph{tilde} on the dimensionless variables for brevity and all subsequent results
 are to be understood
 in terms of the dimensionless variables.
 From (\ref{eq:7}) and (\ref{eq:7new}), the admissible $\Qvec$-tensors belong to the space
 \begin{eqnarray} && \Acal_{\Qvec} = \left\{\Qvec\in
W^{1,2}\left(B(0,R_t), S_0\right); \textrm{$\Qvec=\sqrt{\frac{3}{2}}
 \left(\frac{\xvec\otimes \xvec}{|\xvec|^2}  -
\frac{1}{3}\mathbf{I}\right)$ on $\partial B(0,R_t)$}
\right\}\label{eq:nondim5}
\end{eqnarray}
The associated Euler-Lagrange equations are \cite{maj1,amaz} -
\begin{equation}
\Delta \Qvec_{ij}=  - \Qvec_{ij}  - \frac{ 3\sqrt{6}
h_+}{t}\left(\Qvec_{ik}\Qvec_{kj}-\frac{\delta_{ij}}{3}\textrm{tr}(\Qvec^2)\right)
+\frac{2 h_+^2}{t}\Qvec_{ij}\textrm{tr}(\Qvec^2),~ ~i,j=1,2,3
\label{eq:nondim6}.
\end{equation} All
global and local energy minimizers in $\Acal_\Qvec$ are classical
solutions of (\ref{eq:nondim6}).

\begin{proposition}
\label{prop:u1} For each $t>0$, assume that a  uniaxial global
minimizer of $\tilde{\mathcal{I}}_{LG}$ exists in the admissible space
$\Acal_\Qvec$ defined in (\ref{eq:nondim5}); we denote this uniaxial minimizer by $\Qvec^t$. Let $\left\{t^k\right\}$
be a sequence such that $t^k \to \infty$ as $k\to\infty$ with
$\Qvec^{t_k} \to \Qvec^0$ in $W^{1,2}\left(B(0,R_{t_k}),S_0\right)$ as
$k\to\infty$, where $\Qvec^0$ is the limiting harmonic map
$$ \Qvec^0 = \sqrt{\frac{3}{2}}\left(\frac{\xvec \otimes \xvec}{|\xvec|^2} - \frac{\mathbf{I}}{3} \right). $$ Then from Proposition~\ref{prop:4}, 
\begin{equation}
\lim_{t_k\to \infty} \int_{B(0, R_{t_k})} - \frac{1}{2}\textrm{tr}\tilde{\Qvec}^2 -
\frac{\sqrt{6}h_+}{t_k}\textrm{tr}\tilde{\Qvec}^3 + \frac{h_+^2}{2t_k}
\left(\textrm{tr}\tilde{\Qvec}^2\right)^2 + C(t_k) = 0
\end{equation}
where $R_{t_k} \propto \sqrt{t_k}$. For $k$ sufficiently large, we have
\newline (i) $\Qvec^{t_k} =
\sqrt{\frac{3}{2}}|\Qvec^{t_k}|\left(\nvec\otimes\nvec -
\frac{1}{3}\mathbf{I}\right)$ for some $\nvec \in
W^{1,2}(B(0,R_{t_k}),S^2)$.
\newline (ii) $\Qvec^{t_k}$ is a classical solution of the following
nonlinear system of elliptic partial differential equations on
$B(0,R_{t_k})$ where $R_{t_k}\to \infty$ as $t_k\to\infty$ :
\begin{equation}
\label{eq:24} \Delta\Qvec^{t_k}_{ij} = \left(|\Qvec^{t_k}|^2 -
1\right)\Qvec^{t_k}_{ij} + \frac{3 h_+}{t_k}\left(|\Qvec^{t_k}|^2 -|\Qvec^{t_k}|
\right)\Qvec^{t_k}_{ij}
\end{equation} \newline (iii) $|\Qvec^{t_k}|(\xvec) \leq 1 $ for all $\xvec \in
B(0,R_{t_k})$, \newline (iv) $\frac{1}{R_{t_k}}
\tilde{\mathcal{I}}_{LG}[\Qvec^{t_k}] \leq 12 \pi$ and
\newline (v) $\lim_{t_k \to\infty}\Qvec^{t_k}(0) = 0. $
\newline (vi) All derivatives of $\Qvec^{t_k}$ can be bounded independently of $t_k$ in the $k\to\infty$ limit i.e.
\begin{equation}
\label{eq:gradientbound}
\| \grad^j \Qvec^{t_k} \|_{L^{\infty}\left(B(0,R_{t_k})\right)} \leq C_j \quad j\geq 1
\end{equation}
for a positive constant $C_j$ independent of $t_k$.
\end{proposition}

\textit{Proof:} In what follows, we drop the subscript $k$ for brevity and work in the $t\to\infty$ limit. 
\textit{Proof of (i):}  This follows from the uniaxial character
of $\Qvec^t$ i.e. it has two equal eigenvalues
$$ \Qvec^t(\xvec) = \la(\xvec) \left(\evec(\xvec)\otimes\evec(\xvec) + \fvec(\xvec)\otimes\fvec(\xvec)\right) -
2\la \nvec(\xvec)\otimes \nvec(\xvec) $$ where $\evec,\fvec,\nvec$
form an orthonormal frame at $\xvec$ and $\nvec\otimes\nvec +
\evec\otimes\evec+\fvec\otimes\fvec = \mathbf{I}$. Using the above, $\Qvec^t$ can be
written in the simpler form
\begin{equation}
\label{eq:25} \Qvec^t(\xvec) = -3\la (\xvec)
\left(\nvec(\xvec)\otimes \nvec(\xvec) -
\frac{1}{3}\mathbf{I}\right) \end{equation} where
\begin{equation}
\label{eq:26} |\Qvec^t|^2(\xvec) = 6 \la^2(\xvec).
\end{equation}
From \cite{ejam}, a global uniaxial Landau-de Gennes minimizer has
$\la < 0$ and from (\ref{eq:25})-(\ref{eq:26}), \newline $-3\la =
\sqrt{\frac{3}{2}} |\Qvec^t|$. The representation formula in (i)
follows. $\Box$
\newline \textit{Proof of (ii):} If $\Qvec^t$ is a uniaxial global Landau-de Gennes minimizer,
then it is a classical solution of (\ref{eq:nondim6}). The partial
differential equations (\ref{eq:24}) follow from substituting the
representation formula in (i) into (\ref{eq:nondim6}). $\Box$
\newline \textit{Proof of (iii):} The proof follows from multiplying both sides of (\ref{eq:24})
by $\Qvec^t_{ij}$ and applying a standard maximum principle
argument for $|\Qvec^t|^2$; the details
are omitted for brevity. $\Box$
\newline \textit{Proof of (iv):} This is a direct consequence of
the global energy minimality. The limiting harmonic map $\Qvec^0$
is simply given by
\begin{equation}
\label{eq:27} \Qvec^0 = \sqrt{\frac{3}{2}}\left( \frac{\xvec
\otimes \xvec}{|\xvec|^2} - \frac{1}{3}\mathbf{I} \right)
\end{equation}
in terms of the dimensionless variables in (\ref{eq:nondim}). A
direct computation shows that $\tilde{\mathcal{I}}_{LG}[\Qvec^0] =
12 \pi R_t $ (since $-
\frac{1}{2}\textrm{tr}\left(\Qvec^0\right)^2 -
\frac{\sqrt{6}h_+}{t}\textrm{tr}\left(\Qvec^0\right)^3 +
\frac{h_+^2}{2t} \left(\textrm{tr}\left(\Qvec^0\right)^2\right)^2
+ C(t)=0$ by definition because $\Qvec^0 \in \Qvec_{\min}$) and
hence
\begin{equation}
\label{eq:28} \tilde{\mathcal{I}}_{LG}[\Qvec^t] \leq
\tilde{\mathcal{I}}_{LG}[\Qvec^0] = 12 \pi R_t.
\end{equation} The inequality in (iv) follows. $\Box$
\newline \textit{Proof of (v):} This follows from
Proposition~\ref{prop:4}. We have a topologically non-trivial
boundary condition $\Qvec_b$ in (\ref{eq:8}) and hence every
interior extension of $\frac{\xvec}{|\xvec|}$ must have interior
discontinuities. The extension $\nvec$ in (i) has interior
discontinuities and at every such point of discontinuity
$\xvec^*$, $\Qvec^t(\xvec^*) = 0$ (see \cite{am2} for further
discussion on these lines; $\Qvec^t$ is analytic at $\xvec^*$
whereas $\nvec$ is not and $\nvec$ can lose regularity only when the
number of distinct eigenvalues of $\Qvec^t$ changes. The number of
distinct eigenvalues of $\Qvec^t$ can change only when $\Qvec^t$
relaxes into the isotropic phase i.e. $\Qvec^t(\xvec^*)=0$.) From
Proposition~\ref{prop:4}, as $t\to\infty$, all isotropic points are concentrated
near the singular set of $\Qvec^0$ and the singular set of
$\Qvec^0$ merely consists of the origin. Hence, we have $\lim_{t\to\infty} \Qvec^t(0)=0$. $\Box$

\begin{lemma}
\label{lem:bulk}Assume that  $\Qvec^t$ is a uniaxial global Landau-de
Gennes minimizer in the admissible space $\Acal_\Qvec$ on the
droplet $B(0,R_t)$ in the limit $t\to\infty$. Then the bulk energy
density satisfies the following inequality -
\begin{equation}
\label{eq:bulk} f(\Qvec^t) = - 
\frac{1}{2}\textrm{tr}\left(\Qvec^t \right)^2 -
\frac{\sqrt{6}h_+}{t}\textrm{tr}\left(\Qvec^t \right)^3 +
\frac{h_+^2}{2t} \left(\textrm{tr}\left(\Qvec^t \right)^2\right)^2
+ \frac{1}{2} + \frac{h_+}{t} - \frac{h_+^2}{2t} \geq
\frac{\left(1 - |\Qvec^t|^2 \right)^2}{4}.
\end{equation}
\end{lemma}
\textit{Proof:} From Proposition~\ref{prop:u1}, we have $0\leq
|\Qvec^t| \leq 1$ and one can check from the representation
formula in Proposition~\ref{prop:u1} (i) that
$$\textrm{tr}\left(\Qvec^t \right)^3 =
\frac{|\Qvec^t|^3}{\sqrt{6}}.$$ Substitute the above into the definition of $f(\Qvec^t)$ in (\ref{eq:bulk})
\begin{equation}
\label{eq:bulk1} f(\Qvec^t) = -\frac{1}{2}|\Qvec^t|^2 -
\frac{h_+}{t}|\Qvec^t|^3 + \frac{h_+^2}{2t}|\Qvec^t|^4 +
\frac{1}{2} + \frac{h_+}{t} - \frac{h_+^2}{2t}
\end{equation}
and one can check from (\ref{eq:bulk1}) that
$$f(\Qvec^t) \geq \frac{\left( 1 - |\Qvec^t|^2 \right)^2}{4}.$$
$\Box$

Our second main result concerns the characterization of uniaxial
global Landau-de Gennes minimizers if they exist.
\begin{theorem}
\label{thm:2} Assume that $\Qvec^t\in \Acal_\Qvec$ is a uniaxial global
Landau-de Gennes minimizer on the droplet $B(0,R_t)$ in the limit
$t\to\infty$, where $R_t \propto \sqrt{t} \to \infty$ as $t \to
\infty$. Then $\Qvec^t$ is an entire solution of
\begin{equation}
\label{eq:29} \Delta \Qvec^t_{ij} = \left( |\Qvec^t|^2 - 1\right)
\Qvec^t_{ij} + \frac{3 h_+}{t}\left(|\Qvec^t|^2 - |\Qvec^t|
\right) \Qvec^t_{ij}
\end{equation}
with $\Qvec^t(0)=0$ and $\lim_{t\to\infty}
\frac{1}{R_t}\tilde{\mathcal{I}}_{LG}[\Qvec^t] \leq 12\pi$. There
exists a $\mathbf{T}\in \mathcal{O}(3)$ such that
\begin{equation}
\label{eq:30} \Qvec^t(\xvec) = h(|\xvec|) \left[ \frac{
\mathbf{T}\xvec \otimes \mathbf{T}\xvec}{|\xvec|^2} -
\frac{1}{3}\mathbf{I} \right]
\end{equation}
where $h:\left[0, \infty\right) \to \Rr^+$ is the unique,
monotonically increasing solution of the ordinary differential
equation
\begin{equation}
\label{eq:31} \frac{d^2 h}{d r^2} + \frac{2}{r}\frac{dh}{dr} -
\frac{6 h}{r^2} = h^3 - h + \frac{3 h_+}{t}\left( h^3 - h^2
\right)
\end{equation}
(with $r=|\xvec|$) subject to the boundary conditions
\begin{eqnarray}
\label{eq:32} && h(0)=0 \quad h(r) \to 1 \quad r\to\infty
\nonumber \\ && 0 \leq h(r) \leq 1
\end{eqnarray}
in the limit $t\to\infty$.
\end{theorem}

The last ingredient of the proof of Theorem~\ref{thm:1} is the
following result from \cite{mkaddem&gartland2,am2}:
\begin{theorem}
\label{thm:3} \cite{mkaddem&gartland2,am2}: The radial-hedgehog
solution
\begin{equation}
\label{eq:33} \Hvec = h(|\xvec|) \left[ \frac{ \xvec \otimes
\xvec}{|\xvec|^2} - \frac{1}{3}\mathbf{I} \right]
\end{equation}
where $h$ is the unique solution of (\ref{eq:31}) subject to
(\ref{eq:32}), is a solution of the Euler-Lagrange equations
(\ref{eq:24}). Moreover, $\Hvec$ is not a global Landau-de Gennes
minimizer in the admissible space $\Acal_\Qvec$ in the limit
$t\to\infty$.
\end{theorem}
\textit{Proof:} One can explicitly construct a biaxial perturbation of the form
$$ \mathbf{H}_b(\xvec)  = \Hvec(\xvec) + \frac{1}{(r^2 + 12 )^2}\left(1 - \frac{r}{\sigma} \right) \left( \mathbf{z}\otimes \mathbf{z} - \frac{\mathbf{I}}{3} \right)
$$
with $\sigma=10$ and $\mathbf{z} = (0,0,1)$, that
has lower Landau-de Gennes energy than the radial-hedgehog
solution in the limit $t\to\infty, R_t\to\infty$. This perturbation does not have an isotropic core at the origin. The details can
be found in \cite{am2,mkaddem&gartland2}. $\Box$

To summarize, we obtain symmetry results for global uniaxial
Landau-de Gennes minimizers in the limit $t\to\infty$ i.e. we show
that every uniaxial Landau-de Gennes minimizer is the
radial-hedgehog solution in (\ref{eq:33}) (modulo a rotation) in
the limit $t\to\infty$. The final step is to use results on the
radial-hedgehog solution in the limit $t\to\infty$ from
\cite{mkaddem&gartland2,am2} to prove Theorem~\ref{thm:1}. In the
subsequent sections, we proceed with the proof of
Theorem~\ref{thm:2}.

\section{Symmetry of uniaxial Landau-de Gennes minimizers}
\label{sec:symmetry}

The main trick is to adapt the division trick of \cite{pisante} to
the Landau-de Gennes framework for nematic liquid crystals. Let
$\Qvec^t$ denote a global uniaxial Landau-de Gennes minimizer in
the space $\Acal_\Qvec$ for $t$ sufficiently large so that
Proposition~\ref{prop:u1} holds. Define
\begin{equation}
\label{eq:S} \Svec_{ij}(\xvec) = \frac{\Qvec^t_{ij}(\xvec)}{h(|\xvec|)}
\end{equation}
where $h$ is the unique solution of (\ref{eq:31}) subject to
(\ref{eq:32}) in the $t\to\infty$ limit. We first prove some
auxiliary properties of the function $h:\left[0,\infty\right) \to
\left[0, 1 \right)$ in (\ref{eq:31}).
\begin{theorem}
\label{thm:4} Define the radial-hedgehog solution $\Hvec$ as in
(\ref{eq:33}); $\Hvec$ is a solution of the system of partial
differential equations (\ref{eq:24}) for all $t>1$. In the limit
$t\to\infty$,
\newline (i) $h:\left[0,\infty \right) \to \left[0, 1 \right)$ is
analytic;
\begin{equation}
\label{eq:34} h(0) = h^{'}(0) = 0; h^{''}(0) > 0;
\end{equation} $\Hvec_{ij}(0) = 0$ and $\grad \Hvec_{ij}(0) = 0$,
\newline (ii) $0\leq h(r) \leq 1$ and
\begin{equation}
\label{eq:35} \frac{r^2}{r^2 + 14} \leq h(r) \leq \frac{r^2}{ r^2
+ t \la_t^2}
\end{equation}
where $t \la_t^2 \leq \frac{3}{t}$,
\newline (iii) we have the following gradient bound
\begin{equation}
\label{eq:36} \left| h^{'} (|\xvec|) \right| \leq C_1 \| \grad
\Hvec \|_{L^{\infty}(\Rr^3)} \leq C_2
\end{equation}
where $C_1$ and $C_2$ are positive constants independent of $t$ as
$t\to\infty$,
\newline (iv) we have the following bound for the second
derivative
\begin{equation}
\label{eq:37}\left| h^{''} (|\xvec|) \right| \leq C_3 \| \grad^2
\Hvec \|_{L^{\infty}(\Rr^3)} \leq C_4
\end{equation} where $C_3$ and $C_4$ are positive constants independent of $t$ as
$t\to\infty$ and
\newline (v) we have the following bound for the third derivative
(noting that $\left| \frac{\Hvec_{ij}}{h}\right| = 1$)
\begin{equation}\label{eq:38}
\left| h^{'''}(|\xvec|) \right| =
\left|\frac{\Hvec_{ij}\Hvec_{ij,\alpha\beta\gamma}}{h}
\frac{\xvec_\alpha \xvec_\beta \xvec_\ga}{r^3} \right| \leq C_5\|
\grad^3 \Hvec \|_{L^{\infty}(\Rr^3)} \leq C_6
\end{equation} where $C_5$ and $C_6$ are positive constants independent of $t$ as
$t\to\infty$.
\end{theorem}

\textit{Proof of (i):} The analyticity of $h$, the relations
(\ref{eq:34}) and $\Hvec_{ij}(0)=0$ have been proven in
\cite{am2}. To prove $\grad \Hvec_{ij}(0)=0$, we use the following
equality
\begin{equation}
\label{eq:39} \Hvec_{ij}(\xvec)\Hvec_{ij}(\xvec) = h^2(|\xvec|),~
i,j,=1\ldots 3
\end{equation} so that for any fixed direction characterized by
the unit-vector $\evec_\alpha$ where $\alpha=1\ldots 3$, we have
$$ \Hvec_{ij}(\xvec) \Hvec_{ij,\alpha}(\xvec) = h(|\xvec|)
h^{'}(|\xvec|)\frac{\xvec_\alpha}{|\xvec|}$$ where $\xvec_\alpha =
\xvec \cdot \evec_\alpha$. We set $\xvec = |\xvec|\evec_\alpha$
and re-write the above as
\begin{equation}
\label{eq:40} \lim_{|\xvec|\to 0}
\frac{\Hvec_{ij}\left(|\xvec|\evec_\alpha \right) -
\Hvec_{ij}(0)}{|\xvec|} \Hvec_{ij,\alpha}(|\xvec|\evec_\alpha) =
\lim_{|\xvec|\to 0}\left[ \frac{h(|\xvec|) - h(0)}{|\xvec|}
\right] h^{'}(|\xvec|)
\end{equation}
which implies $\Hvec_{ij,\alpha}(0)\Hvec_{ij,\alpha}(0) = \left[
h^{'}(0)\right]^2 = 0$ for each $\alpha=1\ldots 3$. Hence, $\grad
\Hvec_{ij}(0) = 0$. $\Box$
\newline \textit{Proof of (ii):} The proof of (ii) can be found in
\cite{mkaddem&gartland2} and \cite{am2}.
\newline \textit{Proof of (iii):} We start with the relation
(\ref{eq:39}) and differentiate both sides with respect to the
coordinate direction $\evec_\alpha$ to find
\begin{equation}
\label{eq:41} h^{'}(|\xvec|) =
\frac{\Hvec_{ij}\Hvec_{ij,\alpha}}{h}\frac{\xvec_\alpha}{|\xvec|}.
\end{equation}
Here $\left|\frac{\xvec_\alpha}{|\xvec|}\right| \leq 1$ and
$\left|\frac{\Hvec_{ij}}{h}\right| = 1$. From (\ref{eq:41}), we
have that
\begin{equation}
\label{eq:42} \left| h^{'}(|\xvec|) \right| \leq C_1 || \grad
\Hvec ||_{L^{\infty}(\Rr^3)}
\end{equation}
where $C_1>0$ is independent of $t$. From the gradient bound
derived in \textit{Lemma A.1} in \cite{bbh} and standard results
in elliptic regularity (also refer to Proposition~\ref{prop:u1} (vi)), we have that $|| \grad \Hvec
||_{L^{\infty}(\Rr^3)}$ and all higher derivatives of the
radial-hedgehog solution $\Hvec$ can be bounded independently of
$t$ in the $t\to\infty$ limit. The inequality (\ref{eq:36}) now
follows. $\Box$
\newline \textit{Proof of (iv):}
We start with the relation (\ref{eq:39}) and differentiate twice
to obtain
\begin{equation}
\label{eq:43} \Hvec_{ij}\Hvec_{ij,\alpha\beta} +
\Hvec_{ij,\alpha}\Hvec_{ij,\beta} = \frac{\partial h}{\partial
\xvec_\alpha}\frac{\partial h}{\partial \xvec_\beta} + h
\frac{\partial^2 h}{\partial \xvec_\alpha \partial \xvec_\beta}.
\end{equation}
A direct computation (see (\ref{eq:33})) shows that
$\Hvec_{ij,\alpha}\Hvec_{ij,\beta} =\frac{\partial h}{\partial
\xvec_\alpha}\frac{\partial h}{\partial \xvec_\beta}$ so that
\begin{equation}
\label{eq:44} \frac{\partial^2 h}{\partial \xvec_\alpha \partial
\xvec_\beta} = \frac{\Hvec_{ij}\Hvec_{ij,\alpha\beta}}{h}
\end{equation}
and hence
$$\left| \frac{\partial^2 h}{\partial \xvec_\alpha \partial
\xvec_\beta} \right| \leq C_3 \| \grad^2 \Hvec
\|_{L^{\infty}(\Rr^3)} \leq C_4 $$ for positive constants $C_3$
and $C_4$ independent of $t$ in the $t\to\infty$ limit. Finally,
it suffices to note that
$$ h^{''}(|\xvec|) = \frac{\partial^2 h}{\partial \xvec_\alpha \partial
\xvec_\beta} \frac{\xvec_\alpha \xvec_\beta}{|\xvec|^2}$$ and
(\ref{eq:37}) now follows. $\Box$
\newline \textit{Proof of (v):} We compute an explicit expression
for $\frac{\partial^3 h}{\partial \xvec_\alpha \partial\xvec_\beta
\partial \xvec_\gamma}$ as shown below -
\begin{eqnarray}
\label{eq:45} && \frac{\partial^3 h}{\partial \xvec_\alpha
\partial\xvec_\beta
\partial \xvec_\gamma} = h^{'''}(|\xvec|)\frac{\xvec_\alpha
\xvec_\beta\xvec_\gamma}{|\xvec|^3} +
h^{''}(|\xvec|)\left[\frac{\xvec_\beta \delta_{\alpha\ga} +
\xvec_\alpha \delta_{\beta \ga}}{|\xvec|^2} - \frac{2\xvec_\alpha
\xvec_\beta\xvec_\gamma}{|\xvec|^4}\right] + \nonumber \\
&& + h^{''}(|\xvec|)\frac{\xvec_\ga}{|\xvec|}\left[
\frac{\delta_{\alpha \beta}}{|\xvec|} - \frac{\xvec_\alpha
\xvec_\beta}{|\xvec|^3}\right] + h^{'}(|\xvec|) \left[
-\frac{\delta_{\alpha\beta}\xvec_\ga
+\delta_{\alpha\ga}\xvec_\beta + \delta_{\ga
\beta}\xvec_\alpha}{|\xvec|^3} + \frac{3 \xvec_\alpha
\xvec_\beta\xvec_\gamma}{|\xvec|^5} \right].
\end{eqnarray} We multiply both sides of (\ref{eq:45}) by $\frac{\xvec_\alpha
\xvec_\beta\xvec_\gamma}{|\xvec|^3}$ to obtain
\begin{eqnarray}
\label{eq:46} h^{'''}(|\xvec|) = \frac{\partial^3 h}{\partial
\xvec_\alpha
\partial\xvec_\beta
\partial \xvec_\gamma}\frac{\xvec_\alpha
\xvec_\beta\xvec_\gamma}{|\xvec|^3}.
\end{eqnarray} Straightforward computations show that
\begin{equation}
\label{eq:47} \frac{\partial^3 h}{\partial \xvec_\alpha
\partial\xvec_\beta
\partial \xvec_\gamma} =
\frac{\Hvec_{ij}\Hvec_{ij,\alpha\beta\ga}}{h}
\end{equation}
and on combining (\ref{eq:46}) and (\ref{eq:47}), we obtain
\begin{equation}
\label{eq:48} h^{'''}(|\xvec|) =
\frac{\Hvec_{ij}\Hvec_{ij,\alpha\beta\ga}}{h}\frac{\xvec_\alpha
\xvec_\beta\xvec_\gamma}{|\xvec|^3}.
\end{equation}
The inequality (\ref{eq:38}) now follows. $\Box$

\begin{lemma}
\label{lem:1} Assume that $\Qvec^t$ is a global uniaxial Landau-de Gennes
minimizer in the admissible space $\Acal_\Qvec$ defined in
(\ref{eq:nondim5}) with $\Qvec^t(0)=0$. Then $\grad \Qvec^t(0) =
0$.
\end{lemma}
\textit{Proof:} From the results in \cite{ejam}, we have that a
uniaxial Landau-de Gennes minimizer can be written in the form
$$\Qvec^t = \sqrt{\frac{3}{2}}|\Qvec^t|\left(\nvec\otimes\nvec -
\frac{1}{3}\mathbf{I}\right)$$ for some $\nvec \in
W^{1,2}(B(0,R_t),S^2)$. Therefore, $|\Qvec^t|$ has a minimum at
the origin and we have $\grad |\Qvec^t|(0) = 0$. Using the
relation $\Qvec^t_{ij} \Qvec^t_{ij} = |\Qvec^t|^2$ ,
$\Qvec^t(0)=0$ and $\grad |\Qvec^t|(0) = 0$, we can repeat the
same steps as in Theorem~\ref{thm:4} (i) (see (\ref{eq:39}) and
(\ref{eq:40})) to deduce that $\grad \Qvec^t(0) = 0$. $\Box$

\begin{proposition}
\label{prop:symmetry}Assume that $\Qvec^t$ is a global uniaxial
Landau-de Gennes minimizer in the space $\Acal_\Qvec$ for $t$
sufficiently large, so that Proposition~\ref{prop:u1} holds. Define
\begin{equation}
\Svec_{ij}(\xvec) = \frac{\Qvec^t_{ij}(\xvec)}{h(|\xvec|)} \quad
i,=1\ldots 3 \label{eq:49}.
\end{equation} Then
\begin{equation}\label{eq:50}
|\Svec (\xvec)| \leq C_7 \quad \forall \xvec\in B(0,R_t)\end{equation} where $C_7$ is a  positive
constant independent of $t$. We have the following bounds for
$\grad \Svec$ :
\begin{equation}
\label{eq:51} | \grad \Svec (\xvec)| \leq \begin{cases}
\frac{C_9}{|\xvec|} \quad |\xvec|^2 \leq 1
\\  C_{10} \quad |\xvec|^2
> 1
\end{cases} \end{equation} where $C_9$ and $C_{10}$ are positive
constants independent of $t$.
\end{proposition}
\textit{Proof:} From Theorem~\ref{thm:4}, we have \begin{equation}
\label{eq:54} h (|\xvec|) \geq \frac{|\xvec|^2}{15} \quad
|\xvec|^2 \leq 1\end{equation} so that
$$ |\Svec| = \left| \frac{\Qvec^t}{h} \right| \leq \frac{15
|\Qvec^t|}{|\xvec|^2}$$ where $\| \Qvec^t\|_{L^{\infty}(\Rr^3)}
\leq 1$, $\Qvec^t(0) = 0$ and $\grad \Qvec^t(0) = 0$ (see
Lemma~\ref{lem:1}). From standard Taylor expansion formulae
\cite{taylor}, we have
\begin{equation}
\label{eq:52} \Qvec^t_{ij}(\xvec) = \int_{0}^{1}
\Qvec^t_{ij,\alpha\beta}(s \xvec) \xvec_\alpha \xvec_\beta (1 -
s)~ds \quad |\xvec|^2 \leq 1
\end{equation}
so that
\begin{equation}
\label{eq:53} \left| \Qvec^t_{ij}(\xvec) \right| \leq \| \grad^2
\Qvec^t \|_{L^{\infty}(\Rr^3)} |\xvec|^2.
\end{equation}
$\Qvec^t$ is a classical solution of (\ref{eq:24}) and from
\cite{bbh} and standard results in elliptic regularity
\cite{evans} (also see Proposition~\ref{prop:u1} (vi)), we have that $\|\grad \Qvec^t\|_{L^\infty},
\|\grad^2 \Qvec^t\|_{L^\infty}$ and all higher derivatives of
$\Qvec^t$ can be bounded independently of $t$ in the $t\to\infty$
limit. Substituting the bounds (\ref{eq:54}) and (\ref{eq:53})
into the definition of $\Svec$, we obtain
\begin{equation}
\label{eq:55} \left| \Svec (\xvec) \right| \leq C_{11} \quad
|\xvec|^2 \leq 1
\end{equation}
where $C_{11}$ is a positive constant independent of $t$. For
$|\xvec|^2 \geq 1$, it suffices to note that
$\|\Qvec^t\|_{L^{\infty}(\Rr^3)} \leq 1$ and $h(|\xvec|) \geq
\frac{1}{15}$ for $|\xvec|^2 \geq 1$. Hence, $\left| \Svec (\xvec)
\right| \leq C_{12}$ for $|\xvec|^2 \geq 1$ where $C_{12}$ is a
positive constant independent of $t$ in the $t\to \infty$ limit.
The inequality (\ref{eq:50}) now follows.

Next, we compute bounds for $\grad \Svec$. Consider the case
$|\xvec|^2\leq 1$ first. We re-write $\Svec$ as
\begin{equation}
\label{eq:56} \Svec = \frac{\Qvec^t / |\xvec|^2}{ h(|\xvec|)/
|\xvec|^2}
\end{equation}
so that \begin{equation} \label{eq:60} \left| \grad \Svec \right|
\leq D \left[ \frac{\left| \grad \left(\Qvec^t / |\xvec|^2\right)
\right|}{h(|\xvec|)/ |\xvec|^2} + \frac{\left|\Qvec^t / |\xvec|^2
\right|}{\left( h(|\xvec|)/ |\xvec|^2 \right)^2}\left| \grad
\left(h(|\xvec|)/ |\xvec|^2 \right) \right| \right] \end{equation}
where $D>0$ is a constant independent of $t$. We start with the
integral formula (\ref{eq:52}) and compute $\grad \left(\Qvec^t /
|\xvec|^2\right)$. Straightforward but tedious computations show
that
\begin{equation}
\label{eq:57} \left| \grad \frac{\Qvec^t}{|\xvec|^2} \right| \leq
\frac{D_1}{|\xvec|} \quad |\xvec|^2 \leq 1
\end{equation} where $D_1>0$ is independent of $t$.
Similarly, we have for $|\xvec|^2 \leq 1$ (see
Theorem~\ref{thm:4}),
\begin{equation}
\label{eq:58} h(|\xvec|) = \int_{0}^{|\xvec|} h^{''}(s) \left(
|\xvec| - s \right)~ ds
\end{equation}
so that
\begin{equation}
\label{eq:59} \left| \grad \frac{h(|\xvec|)}{|\xvec|^2}\right|
\leq \frac{D_2}{|\xvec|} \quad |\xvec|^2 \leq 1
\end{equation}where $D_2>0$ is independent of $t$. From the bounds
in Theorem~\ref{thm:4} and (\ref{eq:53}), we have that
$\frac{\Qvec^t}{|\xvec|^2}$ and $\frac{h(|\xvec|)}{|\xvec|^2}$ can
be bounded independently of $t$ for $|\xvec|^2\leq 1$.
Substituting (\ref{eq:57}) and (\ref{eq:59}) into (\ref{eq:60}),
we have
$$ |\grad \Svec(\xvec)| \leq \frac{D_3}{|\xvec|} \quad |\xvec|^2 \leq
1.$$ For $|\xvec|^2 \geq 1$, we have
$$ \frac{\partial \Svec_{ij}}{\partial \xvec_\ga} =
\frac{\Qvec^t_{ij,\ga}}{h} - \frac{\Qvec^t_{ij}}{h^2}
h^{'}(|\xvec|)\frac{\xvec_\ga}{|\xvec|}.$$ For $|\xvec|^2 \geq 1$,
$h$ is bounded away from zero, $\| \Qvec^t \|_{L^{\infty}(\Rr^3)}
\leq 1$, $\|\grad \Qvec^t \|_{L^{\infty}(\Rr^3)}$ and $\|
h^{'}(|\xvec|)\|_{L^{\infty}(\Rr^3)}$ can be bounded independently
of $t$. Hence
$$ |\grad \Svec | \leq D_4  \quad |\xvec|^2 \geq 1 $$
for a constant $D_4 > 0$ independent of $t$. The bounds
(\ref{eq:51}) now follow. $\Box$

We need more careful estimates for $\Svec$ and its gradient near
the origin for the subsequent analysis.
\begin{proposition}
\label{prop:estimate} Let $\Qvec^t$ denote a global uniaxial
Landau-de Gennes minimizer in the space $\Acal_\Qvec$ for $t$
sufficiently large, so that Proposition~\ref{prop:u1} holds. Define $\Svec$ as in (\ref{eq:49}). Then
\begin{eqnarray}
\label{eq:61} && \Svec_{ij}(\xvec)  = \Bvec_{ij\alpha\beta}
\frac{\xvec_\alpha \xvec_\beta}{|\xvec|^2} + o(1) \quad |\xvec|
\to 0 \nonumber \\ && \frac{\partial \Svec_{ij}}{\partial
\xvec_\ga} = \frac{\partial}{\partial \xvec_\ga}
\left[\Bvec_{ij\alpha\beta} \frac{\xvec_\alpha
\xvec_\beta}{|\xvec|^2} \right] + O(1) \quad |\xvec|\to 0
\end{eqnarray} where $\Bvec_{ij\alpha\beta} =
\frac{\Qvec^t_{ij,\alpha\beta}(0)}{h^{''}(0)}$,
$i,j,\alpha,\beta=1\ldots 3$ is a constant matrix.
\end{proposition}

\textit{Proof:} The proof of Proposition~\ref{prop:estimate}
follows from lengthy and technical computations. We state the main
points for completeness. An explicit computation shows that
\begin{equation}
\label{eq:62} \frac{\partial \Svec_{ij}}{\partial \xvec_\ga} =
\frac{\partial}{\partial \xvec_\ga} \left[\Bvec_{ij\alpha\beta}
\frac{\xvec_\alpha \xvec_\beta}{|\xvec|^2} \right] -
A^*_{ij\gamma}\left(\Qvec^t;h\right) + B^*_{ij\gamma}\left(\Qvec^t;h\right) +
C^*_{ij\gamma}\left(\Qvec^t;h\right)
\end{equation}
where
\begin{eqnarray}
\label{eq:63} && A^*_{ij\gamma}\left(\Qvec^t;h\right) =
\frac{\Qvec^t_{ij}/|\xvec|^2}{\left(h(|\xvec|)/|\xvec|^2\right)^2}\frac{\partial}{\partial\xvec_\ga}\left(\frac{h(|\xvec|)}{|\xvec|^2}\right)
\nonumber \\ && B^*_{ij\gamma}\left(\Qvec^t;h\right) =
\frac{|\xvec|^2}{h(|\xvec|)}\frac{\partial}{\partial
\xvec_\ga}\left[ \frac{\Qvec^t_{ij}}{|\xvec|^2} -
\Qvec^t_{ij,\alpha\beta}(0)\frac{\xvec_\alpha \xvec_\beta}{2
|\xvec|^2} \right] \nonumber \\
&& C^*_{ij\gamma}\left(\Qvec^t;h\right) = \frac{\partial}{\partial
\xvec_\ga}\left[\Qvec^t_{ij,\alpha\beta}(0)\frac{\xvec_\alpha
\xvec_\beta}{2 |\xvec|^2}\right]\left[\frac{|\xvec|^2}{h(|\xvec|)}
- \frac{2}{h^{''}(0)}\right].
\end{eqnarray} The next step is to show that $A^*, B^*, C^*$ are bounded
as $|\xvec|\to 0$ in the $t\to\infty$ limit.

We recall that
$$ \frac{h}{|\xvec|^2} \geq \frac{1}{15} $$ for $|\xvec|^2 \leq 1$
(see Theorem~\ref{thm:4}) and
$$ \left| \frac{\Qvec^t(\xvec)}{|\xvec|^2} \right| \leq D_5 \|
\grad^2 \Qvec^t \|_{L^{\infty}(\Rr^3)} $$ where $D_5$ is a
constant independent of $t$ (see (\ref{eq:53})). A straightforward
computation shows that
\begin{equation}
\label{eq:64} h(|\xvec|) = h^{''}(0) \frac{|\xvec|^2}{2} +
\int_{0}^{|\xvec|} h^{'''}(\tau) \frac{(|\xvec| -
\tau)^2}{2}~d\tau
\end{equation}
so that
\begin{eqnarray}
\label{eq:65} &&
\frac{\partial}{\partial\xvec_\ga}\left(\frac{h(|\xvec|)}{|\xvec|^2}\right)
=\frac{\partial}{\partial
|\xvec|}\left(\frac{h(|\xvec|)}{|\xvec|^2}\right)
\frac{\xvec_\ga}{|\xvec|} = \nonumber \\
&&  = \frac{\xvec_\ga}{|\xvec|} \int_{0}^{|\xvec|} h^{'''}(\tau)
\left( 1 - \frac{\tau}{|\xvec|}\right)
\frac{\tau}{|\xvec|^2}~d\tau.
\end{eqnarray}
The equality (\ref{eq:65}) implies that
\begin{equation}
\label{eq:66}
\left|\frac{\partial}{\partial\xvec_\ga}\left(\frac{h(|\xvec|)}{|\xvec|^2}\right)\right|
\leq D_6 \| h^{'''}\|_{L^\infty}
\end{equation}
and from (\ref{eq:38}), we have
\begin{equation}
\label{eq:67}
\left|\frac{\partial}{\partial\xvec_\ga}\left(\frac{h(|\xvec|)}{|\xvec|^2}\right)\right|
\leq D_7
\end{equation}
where $D_6, D_7$ are positive constants independent of $t$. The
inequality (\ref{eq:67}) immediately implies that
\begin{equation}
\label{eq:A} \left| A^*_{ij\gamma}\left(\Qvec^t;h\right) \right| \leq D_8
\end{equation}
where $D_8>0$ is independent of $t$.

Next, we turn to $B^*_{ij\gamma}\left(\Qvec^t;h\right)$. Recalling that
$\Qvec^t_{ij}(0) = \grad \Qvec^t_{ij}(0) = 0$, we have
\begin{equation}
\label{eq:69} \frac{\Qvec^t_{ij}(\xvec)}{|\xvec|^2} =
\frac{\Qvec^t_{ij,\alpha\beta}(0)\xvec_\alpha\xvec_\beta}{2
|\xvec|^2} + \int_{0}^{1} \Qvec^t_{ij,\alpha \beta \sigma}(s\xvec)
\frac{ \xvec_\alpha \xvec_\beta \xvec_\sigma}{2|\xvec|^2} (1 -
s)^2~ ds.
\end{equation}
Therefore,
\begin{equation}
\label{eq:70}
\frac{h(|\xvec|)}{|\xvec|^2}B^*_{ij\gamma}\left(\Qvec^t;h\right) =
\int_{0}^{1} (1-s)^2 \left[\Qvec^t_{ij,\alpha\beta \sigma
\gamma}(s\xvec) s  \frac{ \xvec_\alpha \xvec_\beta
\xvec_\sigma}{2|\xvec|^2} + \Qvec^t_{ij,\alpha\beta\sigma}(s\xvec)
\frac{\partial}{\partial \xvec_\ga}\left(\frac{ \xvec_\alpha
\xvec_\beta \xvec_\sigma}{2|\xvec|^2}\right) \right]~ds.
\end{equation} Hence,
\begin{equation}
\label{eq:B} \left|B^*_{ij\gamma}\left(\Qvec^t;h\right) \right| \leq D_9
\left[ \| \grad^4 \Qvec^t \|_{L^\infty(\Rr^3)} |\xvec| + 1 \right]
\leq D_{10}
\end{equation}
where $D_9$ and $D_{10}$ are positive constants independent of
$t$.

Finally, we turn to $C^*_{ij\gamma}\left(\Qvec^t;h\right)$. From
(\ref{eq:64}), we obtain the following inequality
\begin{equation}
\label{eq:72} \left| |\xvec|^2 h^{''}(0) - 2 h(|\xvec|) \right|
\leq D_{11} \| h^{'''}\|_{L^\infty(\Rr^3)} |\xvec|^3
\end{equation}
which when combined with the inequalities
$$ \left| \frac{\partial}{\partial \xvec_\ga}\left(\frac{ \xvec_\alpha
\xvec_\beta}{|\xvec|^2}\right)\right| \leq \frac{D_{12}}{|\xvec|}
\quad \textrm{and}~ h(|\xvec|) \geq \frac{|\xvec|^2}{15} ~ \textrm{for $|\xvec|^2 \leq 1$}
$$ yields
\begin{equation}
\label{eq:C} \left|C^*_{ij\gamma}\left(\Qvec^t;h\right) \right| \leq D_{13}
\end{equation}
where $D_{11},D_{12}$ and $D_{13}$ are positive constants independent of
$t$ \begin{footnote}{From Theorem~\ref{thm:4}, $h^{''}(0)$ is a
constant independent of $t$ for $t$ sufficiently
large.}\end{footnote}.

From (\ref{eq:A}), (\ref{eq:B}) and (\ref{eq:C}), we deduce that
\begin{equation}
\label{eq:D} \frac{\partial \Svec_{ij}}{\partial \xvec_\ga} =
\frac{\partial}{\partial \xvec_\ga} \left[\Bvec_{ij\alpha\beta}
\frac{\xvec_\alpha \xvec_\beta}{|\xvec|^2} \right] + O(1)~ \textrm{as $|\xvec|\to 0$ and $t\to\infty$.}
\end{equation}

Next, we consider an explicit Taylor expansion for $\Svec_{ij}$
near the origin:
\begin{equation}
\label{eq:73} \Svec_{ij}(\xvec) =
B_{ij\alpha\beta}\frac{\xvec_\alpha \xvec_\beta}{|\xvec|^2} +
B_{ij\alpha\beta}\frac{\xvec_\alpha \xvec_\beta}{2 h
(|\xvec|)|\xvec|^2}\left[ |\xvec|^2 h^{''}(0) - 2 h(|\xvec|) \right] +
\frac{1}{2h}\int_{0}^{1} (1 - s)^2 \Qvec^t_{ij,\alpha\beta
\ga}\left(s \xvec \right)\xvec_\alpha \xvec_\beta \xvec_\ga ~ds.
\end{equation}
From (\ref{eq:72}), we obtain
$$ \left| B_{ij\alpha\beta}\frac{\xvec_\alpha \xvec_\beta}{2 h
(|\xvec|)|\xvec|^2}\left[ h^{''}(0) - 2 h(|\xvec|) \right] \right|
\leq D_{14} \| \grad^2 \Qvec^t \|_{L^\infty(\Rr^3)} |\xvec| \leq
D_{15} |\xvec|$$ and
$$\left|\frac{1}{2h}\int_{0}^{1} (1 - s)^2 \Qvec^t_{ij,\alpha\beta
\ga}\xvec_\alpha \xvec_\beta \xvec_\ga ~ds \right| \leq D_{16}\|
\grad^3 \Qvec^t \|_{L^\infty(\Rr^3)} |\xvec| \leq D_{17}|\xvec| $$ where
$D_{14} - D_{17}$ are positive constants independent of $t$.
Hence,
\begin{equation}
\left| \Svec_{ij}(\xvec) -B_{ij\alpha\beta}\frac{\xvec_\alpha
\xvec_\beta}{|\xvec|^2} \right| \leq D_{18} |\xvec| \quad
|\xvec|\to 0
\end{equation}
for a positive constant $D_{18}$ independent of $t$.
Proposition~\ref{prop:estimate} now follows. $\Box$

\begin{theorem}
\label{thm:5} Assume that there exists a global uniaxial Landau-de Gennes minimizer $\Qvec^t \in \Acal_\Qvec$ for each $t>0$. Let
$\left\{t_k\right\}$ be a sequence such that $t_k \to\infty$ as
$k\to\infty$ so that (up to a subsequence), $\Qvec^{t_k} \to
\Qvec^0$ in $W^{1,2}\left(B(0,R_{t_k}), S_0\right)$ as $t_k\to\infty$,
where $R_k \propto \sqrt{t_k}$ and $\Qvec^0$ is the unique
limiting harmonic map defined in (\ref{eq:27}). Define
\begin{equation}
\label{eq:Sk} \left(\Svec_k\right)_{ij}(\xvec) =
\frac{\Qvec^{t_k}_{ij}(\xvec)}{h(|\xvec|)} \quad i,j=1\ldots 3
\end{equation} as in (\ref{eq:49}), where $h:\left[0,\infty\right) \to \left[0,1\right)$ has been defined in (\ref{eq:31}) - (\ref{eq:32}). Then
in the limit $k\to\infty$, $\Svec_k \in C^2 \left( \Rr^3 \setminus
\left\{0\right\}; S_0 \right)$,
\begin{eqnarray}\label{eq:74}
&& \lim_{k \to \infty} \frac{1}{R_k}\int_{B(0,R_k)} \frac{1}{2}|\grad \Svec_k|^2 + \frac{\left(1 - |\Svec_k|^2\right)^2}{4}~dV
 \leq 12 \pi \nonumber \\
&& \lim_{k\to\infty} \frac{1}{R_k} \int_{B(0,R_k)}\frac{| 1 -
|\Svec_k|^2|}{|\xvec|^2}~dV = 0.
\end{eqnarray}
\end{theorem}

\textit{Proof:} In what follows, we drop the subscript $k$ for
brevity and work with $t_k$ large enough so that
Propositions~\ref{prop:u1}, \ref{prop:symmetry} and \ref{prop:estimate} hold. From
global energy minimality and the inequality established in
Lemma~\ref{lem:bulk}, we have (see Proposition~\ref{prop:u1})
\begin{equation}\label{eq:75}
\frac{1}{R_t}\int_{B(0,R_t)} \frac{1}{2}|\grad \Qvec^t|^2 +
\frac{(1 - |\Qvec^t|^2)^2}{4}~dV \leq \frac{1}{R_t}\int_{B(0,R_t)}
\frac{1}{2}|\grad \Qvec^t|^2 + f(\Qvec^t)~dV \leq 12 \pi
\end{equation} where $R_t \propto \sqrt{t}$ and $f(\Qvec^t)$ has been
defined in Lemma~\ref{lem:bulk}. From the strong convergence in
Proposition~\ref{prop:1} and the inequality established in
Lemma~\ref{lem:bulk}, we have
\begin{equation}
\label{eq:75a} \lim_{t\to\infty}\frac{1}{R_t}\int_{B(0,R_t)}
\frac{\left(1 - |\Qvec^t|^2 \right)^2}{4}~dV \leq
\lim_{t\to\infty}\frac{1}{R_t}\int_{B(0,R_t)}
f\left(\Qvec^t\right)~dV = 0.
\end{equation}
We note that
\begin{equation}
\label{eq:76} \left| 1 - |\Svec|^2 \right|^2 \leq 2\left( \left| 1
- \frac{1}{h^2} \right|^2 + \frac{\left| 1 - |\Qvec^t|^2
\right|^2}{h^4} \right).
\end{equation}
For $|\xvec|^2 \leq 1$, $|\Svec| \leq G_1$ for a positive constant
$G_1$ independent of $t$; see Proposition~\ref{prop:symmetry}. For
$|\xvec|^2\geq 1$, we recall the bounds in Theorem~\ref{thm:4} to
find
\begin{eqnarray}
\label{eq:77} && \left| 1 - \frac{1}{h^2} \right|^2 \leq
\frac{\delta_1}{|\xvec|^4} \quad |\xvec|^2 \geq 1 \nonumber \\
&&\left| 1 - |\Svec|^2 \right|^2 \leq \delta_2 \left| 1 -
|\Qvec^t|^2 \right|^2 \quad |\xvec|^2 \geq 1
\end{eqnarray}
where $\delta_1, \delta_2>0$ are independent of $t$. Combining the
inequalities (\ref{eq:76}) and (\ref{eq:77}), we obtain
\begin{equation}
\label{eq:78} \lim_{t\to\infty}\frac{1}{R_t}\int_{B(0,R_t)}\left|
1 - |\Svec|^2 \right|^2~dV  \leq \delta_3
\lim_{t\to\infty}\frac{1}{R_t}\int_{B(0,R_t)}\left(1 - |\Qvec^t|^2
\right)^2~ dV = 0
\end{equation} where $\delta_3>0$ is independent of $t$.

Next, we turn to the elastic term, $|\grad \Svec|^2$. For
$|\xvec|^2 \leq r_0^2$ where $r_0>>1$ is a constant independent of
$t$, we use the estimates in Proposition~\ref{prop:symmetry} to
obtain
\begin{equation}
\label{eq:79} \int_{B(0,r_0)} \frac{1}{2}|\grad \Svec|^2~dV \leq
\delta_4 r_0^3
\end{equation}
where $\delta_4>0$ is independent of $t$. On the region
$B(0,R_t)\setminus B(0,r_0)$, an explicit computation shows that
\begin{equation}
\label{eq:80} |\grad \Svec|^2 = \frac{|\grad \Qvec^t|^2}{h^2} +
|\Qvec^t|^2 \left( \frac{h^{'}}{h^2}\right)^2 -
2\Qvec^t_{ij}\Qvec^t_{ij,k}\frac{\xvec_k}{|\xvec|}
\frac{h^{'}}{h^3} \quad i,j,k=1\ldots 3.
\end{equation}
For $|\xvec|\geq r_0$, we recall from \cite{am2} that
\begin{equation} \label{eq:81} \left| h^{'}(|\xvec|) \right| \leq
\frac{\delta_5}{|\xvec|^3}
\end{equation}
for a positive constant $\delta_5$ independent of $t$ and
$$ \frac{|\grad \Svec|^2}{h^2} = |\grad \Qvec^t|^2 (1 + \frac{\delta_6}{|\xvec|^2}) $$ for a positive constant $\delta_6>0$ independent of $t$, in
this regime. Combining (\ref{eq:78}) - (\ref{eq:81}), we deduce
the following chain of inequalities
\begin{eqnarray}
\label{eq:82}
 && \lim_{t\to\infty}\frac{1}{R_t}\int_{B(0,R_t)}\frac{1}{2}|\grad \Svec|^2 + \frac{\left|
1 - |\Svec|^2 \right|^2}{4} ~dV =
\lim_{t\to\infty}\frac{1}{R_t}\int_{B(0,R_t)}\frac{1}{2}|\grad
\Svec|^2~dV \leq \nonumber \\ && \leq \lim_{t\to\infty}
\frac{1}{R_t}\int_{B(0,R_t)} \frac{1}{2}|\grad \Qvec^t|^2~dV \leq
12 \pi
\end{eqnarray} where the last inequality follows from
Proposition~\ref{prop:u1}.

Finally, we turn to the integral $\lim_{t\to\infty}\frac{1}{R_t}
\int_{B(0,R_t)}\frac{| 1 - |\Svec|^2|}{|\xvec|^2}~dV$. Recall that
$|\Svec|$ satisfies the global upper bound $|\Svec|\leq \delta_7 $
for a positive constant $\delta_7$ independent of $t$ (see
Proposition~\ref{prop:symmetry}). Then a direct computation shows
that
\begin{equation}
\label{eq:83} \int_{B(0,1)}\frac{| 1 - |\Svec|^2|}{|\xvec|^2}~dV
\leq \delta_8
\end{equation}
for a positive constant $\delta_8$ independent of $t$. On the
region $B(0,R_t)\setminus B(0,1)$, we use Young's inequality to
deduce
\begin{equation}
\label{eq:84} \frac{| 1 - |\Svec|^2|}{|\xvec|^2} \leq
\frac{1}{2}\left[ \frac{\left(1 - |\Svec|^2\right)^2}{4} +
\frac{4}{|\xvec|^4} \right].
\end{equation}
Combining (\ref{eq:83})-(\ref{eq:84}) and recalling (\ref{eq:78}),
we obtain
\begin{eqnarray}
\label{eq:85} && \lim_{t\to\infty}\frac{1}{R_t}
\int_{B(0,R_t)}\frac{| 1 - |\Svec|^2|}{|\xvec|^2}~dV \leq
\nonumber \\ && \leq \lim_{t \to \infty}
\left[\frac{\delta_8}{R_t} + \frac{1}{2 R_t}
\int_{B(0,R_t)\setminus B(0,1)}\frac{\left(1 - |\Svec|^2
\right)^2}{4}~dV + \frac{1}{2 R_t} \int_{B(0,R_t)\setminus B(0,1)}
\frac{4}{|\xvec|^4}~dV \right] = 0
\end{eqnarray}
as required. The proof of Theorem~\ref{thm:5} is now complete.
$\Box$

\subsection{ PDE-methods for $\Svec$}
\label{sec:pde}
 Let $\Qvec^t\in \Acal_\Qvec$ denote a uniaxial global Landau-de Gennes minimizer for $t$ sufficiently large so that Propositions~\ref{prop:u1}, \ref{prop:symmetry}, \ref{prop:estimate}  and Theorem~\ref{thm:5} hold. Define $\Svec$ as in (\ref{eq:49}). We recall that $\Qvec^t$ is a classical solution of
 \begin{equation}
 \label{eq:Q}
 \Delta \Qvec^t_{ij} = \left[ \left(|\Qvec^t|^2 - 1 \right) + \frac{3h_+}{t}\left( |\Qvec^t|^2 - |\Qvec^t| \right) \right] \Qvec^t_{ij}
 \end{equation}
 and $h$ is a solution of the ordinary differential equation
 \begin{equation}
 \label{eq:h}
 h^{''} + \frac{2}{|\xvec|} h^{'} - \frac{6 h}{|\xvec|^2} = h^3 - h + \frac{3h_+}{t}\left( h^3 - h^2 \right)
 \end{equation} where $h^{'}$ denotes $\frac{d h}{d|\xvec|}$ etc. Straightforward but lengthy manipulations show that $\Svec$ satisfies the following system of partial differential equations - 
 \begin{equation}
\label{eq:86} \Delta \Svec_{ij} + \left(1 + \frac{3 h_+}{t}\right)
h^2 \left(1 - |\Svec|^2 \right) \Svec_{ij} = - 2 \frac{h^{'}}{h}
\Svec_{ij,k} \frac{\xvec_k}{|\xvec|} - \frac{6
\Svec_{ij}}{|\xvec|^2} + \frac{3 h_+}{t} h \Svec_{ij} \left(1 -
|\Svec| \right) \quad i,j=1\ldots 3.
\end{equation}

Following the methods in \cite{pisante}, we multiply both sides of the partial differential equation (\ref{eq:86}) with $\Svec_{ij,k}\frac{\xvec_k}{|\xvec|}$. One can check that
\begin{eqnarray}
\label{eq:87} && \Svec_{ij,k}\frac{\xvec_k}{|\xvec|} \Delta \Svec_{ij} = \frac{1}{|\xvec|} \left( \frac{\partial \Svec_{ij}}{\partial |\xvec|} \right)^2
+ \frac{\partial}{\partial \xvec_p} \left[ -\frac{1}{2}|\grad \Svec|^2 \frac{\xvec_p}{|\xvec|} + \Svec_{ij,k}\frac{\xvec_k}{|\xvec|} \Svec_{ij,p} \right]
\\ &&  \left(1 + \frac{3h_+}{t}\right) h^2(|\xvec|) ( 1 - |\Svec|^2 ) \Svec_{ij}\Svec_{ij,k}\frac{\xvec_k}{|\xvec|} = \nonumber \\
&& \label{eq:88} = \left(1 + \frac{3h_+}{t}\right) \left[ \frac{(1 - |\Svec|^2)^2}{4}\left[ 2 h h^{'} + \frac{ 2h^2}{r} \right]
- \frac{\partial}{\partial \xvec_p}\left( \frac{\xvec_p}{|\xvec|} \frac{h^2 \left(1 - |\Svec|^2\right)^2}{4} \right) \right] \\
&& \label{eq:89} - 2 \frac{h^{'}}{h} \Svec_{ij,k}\frac{\xvec_k}{|\xvec|} \Svec_{ij,p}\frac{\xvec_p}{|\xvec|}  = -  2 \frac{h^{'}}{h} \left( \frac{\partial \Svec}{\partial |\xvec|} \right)^2  ~ \textrm{and} \\ && \label{eq:90} -6 \frac{\Svec_{ij}}{|\xvec|^2} \Svec_{ij,p}\frac{\xvec_p}{|\xvec|}  = \frac{\partial}{\partial \xvec_p} \left[ \frac{3 \xvec_p}{|\xvec|^3} (1 - |\Svec|^2) \right].
\end{eqnarray}
Using (\ref{eq:87}) - (\ref{eq:90}), we obtain
\begin{eqnarray}
&&
\frac{\partial}{\partial \xvec_p} \left[ \frac{1}{2}|\grad \Svec|^2 \frac{\xvec_p}{|\xvec|} - \frac{\Svec_{ij,k} \xvec_k}{|\xvec|} \Svec_{ij,p} + 
\left(1 + \frac{3h_+}{t}\right) \frac{\xvec_p}{|\xvec|} \frac{h^2 \left(1 - |\Svec|^2\right)^2}{4}  + \frac{3 \xvec_p \left(1 - |\Svec|^2 \right)}{|\xvec|^3} \right] =  \nonumber \\ &&
= \frac{1}{|\xvec|} \left( \frac{\partial \Svec_{ij}}{\partial |\xvec|} \right)^2 + \left(1 + \frac{3h_+}{t}\right)  \frac{(1 - |\Svec|^2)^2}{4}\left[ 2 h h^{'} + \frac{ 2h^2}{r} \right] + 2 \frac{h^{'}}{h} \left( \frac{\partial \Svec}{\partial |\xvec|} \right)^2 - \frac{3 h_+}{t} h (1 - |\Svec|) \Svec_{ij} \Svec_{ij,k} \frac{\xvec_k}{|\xvec|}. \label{eq:91}
\end{eqnarray}
We recall from Theorem~\ref{thm:4} and the results in \cite{am2} that $h^{'}>0$ for $|\xvec|>0$, for $t$ sufficiently large.

Define
\begin{equation}
\label{eq:Phi}
\mathbf{\Phi}_p = \frac{1}{2}|\grad \Svec|^2 \frac{\xvec_p}{|\xvec|} - \frac{\Svec_{ij,k} \xvec_k}{|\xvec|} \Svec_{ij,p} + 
\left(1 + \frac{3h_+}{t}\right) \frac{\xvec_p}{|\xvec|} \frac{h^2 \left(1 - |\Svec|^2\right)^2}{4}  + \frac{3\xvec_p (1 - |\Svec|^2)}{|\xvec|^3} \quad p=1\ldots 3.
\end{equation}

\begin{lemma}
\label{lem:interior}
We have
\begin{equation}
\label{eq:92}
\int_{|\xvec|=\delta} \Phi_p \frac{\xvec_p}{|\xvec|} dA \to 12 \pi  \quad \textrm{as $\delta \to 0$}
\end{equation} where $dA$ is the surface area element on the sphere of radius $\delta$ centered at the origin.
\end{lemma}

\textit{Proof:} By the definition of $\Phi_p$ in (\ref{eq:Phi}), we have
\begin{eqnarray}
&&  \int_{|\xvec|=\delta} \Phi_p \frac{\xvec_p}{|\xvec|} dA =  \nonumber \\
&& = \int_{|\xvec|=\delta} \frac{1}{2}|\grad \Svec|^2 - \left(\frac{\partial \Svec}{\partial |\xvec|}\right)^2 + \left(1 + \frac{3h_+}{t}\right) \frac{h^2 (1 - |\Svec|^2)^2}{4} + \frac{ 3 \left(1 - |\Svec|^2 \right)}{|\xvec|^2} dA. \label{eq:93}
\end{eqnarray}
From the estimates in Proposition~\ref{prop:estimate}, we have the following as $\delta\to 0$
\begin{eqnarray}
&& |\grad \Svec|^2 = \left| \grad \left(\frac{\Bvec_{ij\alpha\beta}\xvec_\alpha \xvec_\beta}{|\xvec|^2} \right) \right|^2 + o\left( |\xvec|^{-2} \right) \nonumber \\ &&  \left(\frac{\partial \Svec}{\partial |\xvec|}\right)^2 = o( |\xvec|^{-2}
\nonumber \\ && 1 - |\Svec|^2 = \frac{ |\xvec|^4 - \left|\Bvec_{ij\alpha\beta}\xvec_\alpha \xvec_\beta \right|^2 }{|\xvec|^4} + o(1)
\label{eq:94}
\end{eqnarray}
where $\Bvec_{ij\alpha\beta} = \frac{\Qvec^t_{ij,\alpha \beta}(0)}{h^{''}(0)}$ is a constant matrix and $i,j,\alpha,\beta=1\ldots 3$.

Substituting (\ref{eq:94}) into (\ref{eq:93}), we get
\begin{eqnarray}
 && \int_{|\xvec|=\delta} \Phi_p \frac{\xvec_p}{|\xvec|} dA =  \nonumber \\ && =
 \int_{|\xvec|=\delta} \frac{1}{2} \left| \grad \left(\frac{\Bvec_{ij\alpha\beta}\xvec_\alpha \xvec_\beta}{|\xvec|^2} \right) \right|^2  + 3\left(\frac{ |\xvec|^4 - \left|\Bvec_{ij\alpha\beta}\xvec_\alpha \xvec_\beta \right|^2 }{|\xvec|^6} \right) + o(|\xvec|^{-2})~ dA = \nonumber \\
 && = 12 \pi + o(1) + \int_{|\xvec|=\delta} \frac{1}{2} \left| \grad \left(\frac{\Bvec_{ij\alpha\beta}\xvec_\alpha \xvec_\beta}{|\xvec|^2} \right) \right|^2 - \frac{3}{|\xvec|^2}\left| \frac{\Bvec_{ij\alpha\beta}\xvec_\alpha \xvec_\beta}{|\xvec|^2} \right|^2 ~dA.
 \end{eqnarray}
 A direct computation shows that
 $$\int_{|\xvec|=\delta} \frac{1}{2} \left| \grad \left(\frac{\Bvec_{ij\alpha\beta}\xvec_\alpha \xvec_\beta}{|\xvec|^2} \right) \right|^2 - \frac{3}{|\xvec|^2}\left| \frac{\Bvec_{ij\alpha\beta}\xvec_\alpha \xvec_\beta}{|\xvec|^2} \right|^2 ~dA = 0$$
 for the constant matrix $\Bvec_{ij\alpha\beta}= \frac{\Qvec^t_{ij,\alpha \beta}(0)}{h^{''}(0)}$ and hence, the conclusion of Lemma~\ref{lem:interior} follows. $\Box$

\begin{lemma}
\label{lem:B}
The integral
\begin{equation}
\label{eq:Bnew}\int_{|\xvec|=1} \frac{1}{2} \left| \grad \left(\frac{\Bvec_{ij\alpha\beta}\xvec_\alpha \xvec_\beta}{|\xvec|^2} \right) \right|^2 - \frac{3}{|\xvec|^2}\left| \frac{\Bvec_{ij\alpha\beta}\xvec_\alpha \xvec_\beta}{|\xvec|^2} \right|^2 ~dA = 0
\end{equation}
for any constant $\Bvec_{ij\alpha \beta}$ such that
$\Bvec_{ij\alpha \beta} = \Bvec_{ji\alpha \beta}, ~ \Bvec_{ij\alpha \beta} = \Bvec_{ij \beta \alpha}$ and
$\Bvec_{ij\alpha \alpha}  =\Bvec_{ii \alpha\beta} = 0$.
\end{lemma}

\textit{Proof:} Consider the matrix
\begin{equation}
\label{eq:B2}
\Bvec_{ij\alpha\beta} = \frac{\Qvec^t_{ij,\alpha \beta}(0)}{h^{''}(0)} \quad i,j,\alpha,\beta=1\ldots 3.
\end{equation} Recalling that $\Qvec^t$ is a symmetric, traceless $3\times 3$ matrix which is a classical solution the PDE-system (\ref{eq:Q}) (so that we have $\Qvec^t_{ij,\alpha\beta} = \Qvec^t_{ij,\beta\alpha}$ from equality of mixed partial derivatives), we obtain
\begin{equation}
\label{eq:B3}
\Bvec_{ij\alpha\beta}  =\Bvec_{ji\alpha\beta} ;~ \Bvec_{ij\alpha\beta} = \Bvec_{ij\beta\alpha}; ~ \Bvec_{ij \alpha \alpha} = 0,~ \Bvec_{ii \alpha \beta}=0
\end{equation}
for all $i,j,\alpha,\beta = 1\ldots 3$.

A direct computation shows that
\begin{equation}
\label{eq:B4}
\frac{1}{2} \left| \grad \left(\frac{\Bvec_{ij\alpha\beta}\xvec_\alpha \xvec_\beta}{|\xvec|^2} \right) \right|^2 = \frac{2}{|\xvec|^4} \Bvec_{ijpq}\Bvec_{ijrs}\xvec_q \xvec_s \left( \delta_{rp} - \frac{\xvec_r \xvec_p}{|\xvec|^2} \right)
\end{equation}
and
\begin{equation}
\label{eq:B5} \frac{3}{|\xvec|^2}\left| \frac{\Bvec_{ij\alpha\beta}\xvec_\alpha \xvec_\beta}{|\xvec|^2} \right|^2 = 
\frac{3}{|\xvec|^6} \Bvec_{ijpq} \Bvec_{ijrs} \xvec_p \xvec_q \xvec_r \xvec_s 
\end{equation}
 so that
 \begin{eqnarray}
&& \int_{|\xvec|=1} \frac{1}{2} \left| \grad \left(\frac{\Bvec_{ij\alpha\beta}\xvec_\alpha \xvec_\beta}{|\xvec|^2} \right) \right|^2 - \frac{3}{|\xvec|^2}\left| \frac{\Bvec_{ij\alpha\beta}\xvec_\alpha \xvec_\beta}{|\xvec|^2} \right|^2 ~dA = \nonumber \\
&& = \Bvec_{ijpq} \Bvec_{ijrs} \left[ 2\delta_{rp} \int_{|\xvec| = 1} \xvec_q \xvec_ s dA -  5 \int_{|\xvec|=1} \xvec_p \xvec_q \xvec_r \xvec_s dA \right] \label{eq:B6}
\end{eqnarray} for $i,j,p,q,r,s= 1\ldots 3$.

Using spherical coordinate representation, we can check that
\begin{equation}
\label{eq:B7}
\int_{|\xvec| = 1} \xvec_q \xvec_ s dA = \frac{4\pi}{3}\delta_{qs}
\end{equation}
and
\begin{equation}
\label{eq:B8}
\int_{|\xvec|=1} \xvec_p \xvec_q \xvec_r \xvec_s dA  = \frac{4\pi}{15}\left[ \delta_{pq}\delta_{rs} + \delta_{pr}\delta_{qs} + \delta_{ps}\delta_{qr} \right].
\end{equation}
Substituting (\ref{eq:B7}) and (\ref{eq:B8}) into (\ref{eq:B6}), we obtain
\begin{eqnarray}
&& \int_{|\xvec|=1} \frac{1}{2} \left| \grad \left(\frac{\Bvec_{ij\alpha\beta}\xvec_\alpha \xvec_\beta}{|\xvec|^2} \right) \right|^2 - \frac{3}{|\xvec|^2}\left| \frac{\Bvec_{ij\alpha\beta}\xvec_\alpha \xvec_\beta}{|\xvec|^2} \right|^2 ~dA = \nonumber \\ && = \frac{4\pi}{3}\left[ 2 \Bvec_{ijrs}\Bvec_{ijrs} - \Bvec_{ijpp}\Bvec_{ijss}  - \Bvec_{ijqr}\Bvec_{ijrq} - \Bvec_{ijsr} \Bvec_{ijrs} \right]
\end{eqnarray}
and the right-hand side vanishes by virtue of the properties established in (\ref{eq:B3}). The integral equality (\ref{eq:Bnew}) now follows. $\Box$

\begin{proposition}
\label{prop:S} Let $\Qvec^t\in \Acal_\Qvec$ denote a uniaxial global Landau-de Gennes minimizer for $t$ sufficiently large so that Propositions~\ref{prop:u1}, \ref{prop:symmetry}, \ref{prop:estimate} and Theorem~\ref{thm:5} hold. Define $\Svec$ as in (\ref{eq:49}). Then
\begin{eqnarray}
\label{eq:B9}
&& \frac{\partial \Svec}{\partial |\xvec|} = 0 \quad \xvec \in B(0, R_t )\\
&& |\Svec(\xvec)| = 1 \quad \xvec \in B(0, R_t )
\end{eqnarray}
where $R_t = \mu \sqrt{t}$ and $\mu$ is a constant independent of $t$.
\end{proposition}

\textit{Proof:} We integrate both sides of (\ref{eq:91}) from $|\xvec|=0$ to $|\xvec|=R_t$, divide by $R_t$, use Lemma~\ref{lem:interior} and take limit $t\to\infty$ to obtain
\begin{eqnarray}
&&
12\pi + \lim_{t \to \infty} \frac{1}{R_t}\int_{0}^{R_t}\int_{B(0,R)}\frac{1}{|\xvec|}\left(\frac{\partial \Svec}{\partial |\xvec|}\right)^2
+ \left(1 + \frac{3 h_+}{t} \right) \frac{\left( 1- |\Svec|^2 \right)^2}{4}\left[ 2 h^{'} h + \frac{ 2 h^2}{|\xvec|} \right] + \frac{ 2 h^{'}}{h} \left(\frac{\partial \Svec}{\partial |\xvec|}\right)^2~dV dR - \nonumber \\ && - \lim_{t\to\infty} \frac{3h_+}{t} \frac{1}{R_t}\int_{0}^{R_t}\int_{B(0,R)} h (1 - |\Svec|) \Svec_{ij} \Svec_{ij,k} \frac{\xvec_k}{|\xvec|}~dV~dR = \nonumber \\
&& = \lim_{t \to \infty}  \frac{1}{R_t}\int_{B(0,R_t)} \frac{1}{2}|\grad \Svec|^2 - \left(\frac{\partial \Svec}{\partial |\xvec|}\right)^2
+ \left(1 + \frac{3 h_+}{t} \right) \frac{h^2 \left( 1- |\Svec|^2 \right)^2}{4} + \frac{3\left(1 - |\Svec|^2 \right)}{|\xvec|^2}~dV.
\label{eq:B10}
\end{eqnarray}
From (\ref{eq:74}), we have that 
$$ \lim_{t \to \infty}  \frac{1}{R_t}\int_{B(0,R_t)} \frac{3\left(1 - |\Svec|^2 \right)}{|\xvec|^2}~dV = 0 $$
and
$$ \lim_{t \to \infty}  \frac{1}{R_t}\int_{B(0,R_t)}  \frac{1}{2}|\grad \Svec|^2 - \left(\frac{\partial \Svec}{\partial |\xvec|}\right)^2 + \left(1 + \frac{3 h_+}{t} \right) \frac{h^2 \left( 1- |\Svec|^2 \right)^2}{4} \leq 12 \pi $$
since $h^2(|\xvec|) \leq 1$ on $B(0,R_t)$ and $\lim_{t\to\infty} \frac{h_+}{t} = 0$.

From Theorem~\ref{thm:4}, we recall that $h$ is monotonically increasing in the limit $t\to\infty$ so that every term in the integrand of 
$$ \frac{1}{R_t}\int_{0}^{R_t}\int_{B(0,R)}\frac{1}{|\xvec|}\left(\frac{\partial \Svec}{\partial |\xvec|}\right)^2
+ \left(1 + \frac{3 h_+}{t} \right) \frac{\left( 1- |\Svec|^2 \right)^2}{4}\left[ 2 h^{'} h + \frac{ 2 h^2}{|\xvec|} \right] + \frac{ 2 h^{'}}{h} \left(\frac{\partial \Svec}{\partial |\xvec|}\right)^2~dV dR$$ is non-negative.

Finally, we estimate the integral $\frac{1}{R_t}\int_{0}^{R_t}\int_{B(0,R)} h (1 - |\Svec|) \Svec_{ij} \Svec_{ij,k} \frac{\xvec_k}{|\xvec|}~dV~dR$ as follows:
\begin{eqnarray}
&&
\left| \int_{B(0,R)} h (1 - |\Svec|) \Svec_{ij} \Svec_{ij,k} \frac{\xvec_k}{|\xvec|}~dV \right|
\leq \nonumber \\ && \leq C \left[ \int_{B(0,R)}\left(1 - |\Svec|^2 \right)^2~dV \right]^{1/2} \left[\int_{B(0,R)} |\grad \Svec|^2~dV \right]^{1/2}
\leq  C^* \sqrt{f[R]}\sqrt{R} \label{eq:B11}
\end{eqnarray} where $\lim_{R\to\infty} \frac{f[R]}{R}=0.$
for all $R>0$ where $C$ and $C^*$ are positive constants independent of $R$. The first inequality follows from Cauchy-Schwarz and the fact that $|\Svec|\leq D$ for a positive constant $D$ independent of $R$ (see Proposition~\ref{prop:symmetry}).
From Theorem~\ref{thm:5}, we have
$$ \lim_{R\to\infty} \frac{1}{R} \int_{B(0,R)}\left(1 - |\Svec|^2 \right)^2~dV =0 $$
so that if we set
$$\int_{B(0,R)}\left(1 - |\Svec|^2 \right)^2~dV = f[R],$$ then
$f[R]=o(R)$ as $R\to\infty$. The second inequality also uses the upper bound $\int_{B(0,R)} |\grad \Svec|^2~dV \leq 12\pi R$ from (\ref{eq:74}), yielding the inequality (\ref{eq:B11}). 

Substituting the upper bound (\ref{eq:B11}) into (\ref{eq:B10}), we get
\begin{eqnarray}
\label{eq:B12}
\lim_{t\to\infty} \frac{3h_+}{t} \left| \frac{1}{R_t}\int_{0}^{R_t}\int_{B(0,R)} h (1 - |\Svec|) \Svec_{ij} \Svec_{ij,k} \frac{\xvec_k}{|\xvec|}~dV~dR \right|
 = 0
\end{eqnarray} since $\int_{0}^{R} \sqrt{f[R] R }~dR = o(R^2)$ as $R\to\infty$ and $\frac{3h_+}{t} \sim \frac{1}{R_t}$ as $t\to\infty$.

Combining the above, we deduce that
\begin{equation}
\label{eq:B13} \int_{B(0,R)}\frac{1}{|\xvec|}\left(\frac{\partial \Svec}{\partial |\xvec|}\right)^2
+ \left(1 + \frac{3 h_+}{t} \right) \frac{\left( 1- |\Svec|^2 \right)^2}{4}\left[ 2 h^{'} h + \frac{ 2 h^2}{|\xvec|} \right] + \frac{ 2 h^{'}}{h} \left(\frac{\partial \Svec}{\partial |\xvec|}\right)^2~dV = 0
\end{equation}
for all $R>0$ and the conclusion of Proposition~\ref{prop:S} now follows. $\Box$

\begin{proposition}
\label{prop:main}
Let $\Qvec^t\in \Acal_\Qvec$ denote a uniaxial global Landau-de Gennes minimizer for $t$ sufficiently large so that Propositions~\ref{prop:u1}, \ref{prop:symmetry}, \ref{prop:estimate} and Theorem~\ref{thm:5} hold. Then in the limit $t\to\infty$, we have
\begin{equation}
\label{eq:main}
\Qvec^t = h( |\xvec|) \left( \frac{\xvec \otimes \xvec}{|\xvec|^2} - \frac{\mathbf{I}}{3} \right) \quad \xvec \in B(0,R_t)
\end{equation} where $h$ is the unique solution of (\ref{eq:31}) subject to
(\ref{eq:32}) and $R_t = \mu \sqrt{t}$, with $\mu$ being a constant independent of $t$.
\end{proposition}

\textit{Proof:} From Proposition~\ref{prop:S}, we have that in the limit $t\to \infty$
\begin{equation}
\label{eq:main1}
\Qvec^t = h(|\xvec|) \mathbf{M}_{ij}\left( \frac{\xvec}{|\xvec|} \right)
\end{equation}
where  $\mathbf{M}_{ij} = \sqrt{\frac{3}{2}}\left(\mvec\otimes \mvec - \frac{\mathbf{I}}{3} \right)$ for some $\mvec \in S^2$ (from the uniaxial character of $\Qvec^t$), 
\begin{equation}
\label{eq:main2}
\left| \mathbf{M} (\xvec)\right|^2 = 1 \quad \xvec \in B(0,R_t)
\end{equation}
and $\mathbf{M}_{ij} \to \sqrt{\frac{3}{2}}\left(\frac{\xvec\otimes \xvec}{|\xvec|} - \frac{\mathbf{I}}{3}\right) $ as $|\xvec| \to\infty$ from the imposed Dirichlet boundary condition. 

We substitute (\ref{eq:main1}) into the governing system of PDEs for $\Qvec^t$ as shown below - 
\begin{equation}
\label{eq:main3}
\Delta \Qvec^t_{ij} = \left[ \left( |\Qvec^t|^2 - 1 \right) + \frac{3 h_+}{t}\left( |\Qvec^t|^2 - |\Qvec^t| \right) \right] \Qvec^t_{ij},
\end{equation}
multiply both sides of (\ref{eq:main3}) with $\mathbf{M}_{ij}$ (noting that $\mathbf{M}_{ij} \mathbf{M}_{ij,k}=0$) to find
\begin{equation}
\label{eq:main4}
\left| \grad \mathbf{M} \right|^2 = 3 |\grad \mvec|^2 = \frac{6}{r^2}.
\end{equation}

Consider the minimization problem
\begin{equation}
\label{eq:main5}
\min_{\nvec \in N} \int_{B(0,R)} |\grad \mathbf{n}|^2~dV
\end{equation}
where $B(0,R)$ is a three-dimensional droplet of arbitrary radius $R>0$ and
\begin{equation}
\label{eq:main6}
N = \left\{ \mathbf{n} \in W^{1,2}(B(0,R);S^2); \nvec = \frac{\xvec}{|\xvec|} ~on~\partial B(0,R) \right\}.
\end{equation}
It is known that the minimization problem (\ref{eq:main5})-(\ref{eq:main6}) has a \emph{unique} minimizer \cite{lin}
\begin{equation}
\label{eq:main7}
\nvec_{\min} = \frac{\xvec}{|\xvec|}
\end{equation}
and
\begin{equation}
\label{eq:main8}
\int_{B(0,R)} |\grad \mathbf{n}_{\min}|^2~dV = 8 \pi R.
\end{equation}

The unit-vector field $\mathbf{m} \in  W^{1,2}(B(0,R_t);S^2)$, $\mathbf{m}=\frac{\xvec}{|\xvec|}$ on $\partial B(0,R_t)$ and from (\ref{eq:main4}), we have
\begin{equation}
\label{eq:main9}
\int_{B(0,R)} |\grad \mathbf{m}|^2~dV = 8 \pi R_t.
\end{equation}
Comparing (\ref{eq:main8}) and (\ref{eq:main9}), we deduce that $\mvec$ is a minimizer of the problem (\ref{eq:main5})-(\ref{eq:main6}) on $B(0,R_t)$ and from the uniqueness of the minimizer, we deduce that
$$ \mvec(\xvec) = \frac{\xvec}{|\xvec|} \quad \xvec\in B(0,R_t). $$
The conclusion of Proposition~\ref{prop:main} now follows. $\Box$


\begin{thebibliography}
\smallskip
\par
\bibitem{bz}  J.M. Ball and A. Zarnescu, Orientability and energy
minimization for liquid crystals, in preparation
\bibitem{bbh} F. Bethuel, H. Brezis  and F.H\'{e}lein,
Asymptotics for the minimization of a Ginzburg-Landau functional.
Calc. Var. Partial Differential Equations  1  (1993), no. 2,
123--148
 \bibitem{gartland} T.~Davis and E.~Gartland, Finite element analysis of
the Landau--De Gennes minimization problem for liquid crystals.
SIAM Journal of Numerical Analysis, 35, 336-362 (1998).
 \bibitem{ericksen} J. L. Ericksen,  Liquid crystals with variable degree
of orientation.
 Arch. Rational Mech. Anal.  113  (1990),  no. 2, 97--120
 \bibitem{evans} L.~Evans, Partial Differential Equations. American
Mathematical Society, Providence, 1998.
\bibitem{oseenfrank} F.C. Frank,  On the theory of liquid crystals.
 Disc. Faraday Soc., 25(1958)1
\bibitem{dg} P. G. De Gennes,  The physics of liquid crystals.
Oxford, Clarendon Press. 1974
\bibitem{gilbarg} D.Gilbarg and N.Trudinger, Elliptic Partial Differential
Equations of Second Order.
 Springer, 224, 2, 1977
\bibitem{partialcrystal} R. Hardt, D. Kinderlehrer and F. H. Lin,
Existence and partial regularity of static liquid crystals
configurations, Comm. Math. Phys., 105 (1986), 547-570
 \bibitem{lin} F.~H.~Lin and C.~Liu, Static and Dynamic Theories of Liquid
Crystals. Journal of Partial Differential Equations, 14, no. 4,
289--330  (2001).
\bibitem{ejam}  A. Majumdar, Equilibrium order parameters of liquid crystals
in the Landau--de Gennes theory. European Journal of Applied Mathematics, \textbf{21}, 181-203 (2010).
\bibitem{amaz}A.Majumdar \& A.Zarnescu, The Landau-de Gennes
theory of nematic liquid crystals: the Oseen-Frank limit and
beyond. Archive of Rational Mechanics and
Analysis, \textbf{196}, No 1,   227--280 (2010).
\bibitem{maj1} A.Majumdar, The Landau-de Gennes theory of nematic liquid
crystals: Uniaxiality versus Biaxiality. Under review in Communications in Pure and Applied Analysis.
\bibitem{am2} A. Majumdar, The Radial-Hedgehog Solution in the Landau-de Gennes theory for nematic liquid crystals. Accepted
for publication in the European Journal of Applied Mathematics.
\bibitem{mkaddem&gartland2} S. Mkaddem and E. C. Gartland,
 On the local instability of radial hedgehog configurations in nematic liquid crystals under Landau-de Gennes free-energy models.   Phys.
Rev.Rev. E. \textbf{59} 563--567 (1999).
 \bibitem{moser}  R. Moser,  Partial regularity for harmonic maps and
related problems. World Scientific Publishing , Hackensack, NJ, 2005.
\bibitem{newtonmottram} N.J.Mottram and C.Newton,  Introduction to
\textbf{Q}-tensor Theory.
 University of Strathclyde, Department of Mathematics, Research Report,
10, 2004
\bibitem{pisante} A.Pisante, Two results on the equivariant Ginzburg–Landau vortex in arbitrary dimension. Journal of Functional Analysis , 260 (3), 892 -- 905 (2011). 
\bibitem{schoen} R. Schoen and K. Uhlenbeck,  A Regularity Theory for Harmonic 
Mappings. J. Diff. Geom. 1982, 17, 307-335. 
\bibitem{virga}E. G. Virga, Variational theories for liquid crystals.
Chapman and Hall, London 1994
\end{thebibliography}
\end{document}